\documentclass[10pt]{article}
\usepackage{amsbsy}
\def\Z{\bf \mbox{Z\hspace{-.40em}Z}}
\def\R{{\bf \mbox{I\hspace{-.20em}R}}}

\def\lcf{\lbrack\! \lbrack}
\def\rcf{\rbrack\! \rbrack}

\newtheorem{definition}{Definition}[section]
\newtheorem{lemma}[definition]{Lemma}
\newtheorem{proposition}[definition]{Proposition}
\newtheorem{theorem}[definition]{Theorem}
\newtheorem{remark}[definition]{Remark}
\newtheorem{corol}[definition]{Corollary}
\newtheorem{remarks}[definition]{Remarks}

\begin{document}

\def \biblio{}

\title{REDUCTION OF JACOBI MANIFOLDS VIA DIRAC STRUCTURES THEORY}
\date{}
\maketitle \vspace{-10mm}
\begin{center}
\begin{tabular}{ccc}
{\large Fani PETALIDOU} &  & {\large Joana M. NUNES DA COSTA} \\
{\it Faculty of Sciences and Technology}  & & {\it Department of Mathematics} \\
{\it University of Peloponnese} & & {\it University of Coimbra} \\
{\it 22100 Tripoli, Greece}  & & {\it Apartado 3008} \\
& & {\it 3001-454 Coimbra, Portugal} \\
 &  & \\
{\it e-mail : petalido@uop.gr} & & {\it e-mail :
jmcosta@mat.uc.pt}
\end{tabular}
\end{center}

\vspace{20mm}
\begin{abstract}
We first recall some basic definitions and facts about Jacobi
manifolds, generalized Lie bialgebroids, generalized Courant
algebroids and Dirac structures. We establish an one-one
correspondence between reducible Dirac structures of the
generalized Lie bialgebroid of a Jacobi manifold $(M,\Lambda,E)$
for which $1$ is an admissible function and Jacobi quotient
manifolds of $M$. We study Jacobi reductions from the point of
view of Dirac structures theory and we present some examples and
applications.

\vspace{3mm} \noindent {\bf{Keywords : }}{Dirac structures,
generalized Lie bialgebroids, generalized Courant algebroids,
Jacobi manifolds, reduction.}

\vspace{3mm} \noindent {\bf{A.M.S. classification (2000) :}}
53D10, 53D17, 53D20, 58A30.
\end{abstract}

\section{Introduction}
The concept of a {\it Dirac structure} on a differentiable
manifold $M$ was introduced by T. Courant and A. Weinstein in
\cite{cw} and developed by T. Courant in \cite{c}. Its principal
aim is to present a unified framework for the study of
pre-symplectic forms, Poisson structures and foliations. More
specifically, a {\it Dirac structure} on $M$ is a subbundle
$L\subset TM\oplus T^{\ast}M$ that is maximally isotropic with
respect to the canonical symmetric bilinear form on $TM\oplus
T^{\ast}M$ and satisfies a certain integrability condition. In
order to formulate this integrability condition, T. Courant
defines a bilinear, skew-symmetric, bracket operation on the space
$\Gamma(TM\oplus T^{\ast}M)$ of smooth sections of $TM\oplus
T^{\ast}M$ which does not satisfy the Jacobi identity. The nature
of this bracket was clarified by Z.-J. Liu, A. Weinstein and P. Xu
in \cite{lwx1} by introducing the structure of a {\it Courant
algebroid} on a vector bundle $E$ over $M$ and by extending the
notion of a Dirac structure to the subbundles $L\subset E$. The
most important example of Courant algebroid is the direct sum
$A\oplus A^{\ast}$ of a Lie bialgebroid $(A,A^{\ast})$ over a
smooth manifold $M$ (\cite{mx}).

Alan Weinstein and his collaborators have studied several problems
of Poisson geometry via Dirac structures theory. In \cite{lwx2},
Z.-J. Liu {\it et al.} establish an one-one correspondence between
Dirac subbundles of the double $TM\oplus T^{\ast}M$ of the
triangular Lie bialgebroid $(TM,T^{\ast}M,\Lambda)$ defined by a
Poisson structure $\Lambda$ on $M$ and Poisson structures on
quotient manifolds of $M$. Using this correspondence and the
results concerning the pull-backs Dirac structures under Lie
algebroid morphisms, Z.-J. Liu constructs in \cite{liu} the
Poisson reduction in terms of Dirac structures.

\vspace{1mm}

On the other hand, it is well known that the notion of {\it Jacobi
manifold}, i.e. a differentiable manifold $M$ endowed with a
bivector field $\Lambda$ and a vector field $E$ satisfying an
integrability condition, introduced by A. Lichnerowicz in
\cite{lch}, is a rich geometrical notion that generalizes the
Poisson, symplectic, contact and co-symplectic manifolds. Thus, it
is natural to research a simple interpretation of Jacobi manifolds
by means of Dirac structures. A first approach of this problem is
presented in \cite{wd} by A. Wade. Taking into account that to any
Jacobi structure $(\Lambda,E)$ on $M$ is canonically associated a
generalized Lie bialgebroid structure on $(TM\times \R,
T^{\ast}M\times \R)$ (\cite{im1}), she considers the Whitney sum
$\mathcal{E}^1(M)=(TM\times \R)\oplus (T^{\ast}M\times \R)$,
introduces the notion of $\mathcal{E}^1(M)$-{\it Dirac structures}
by extending the Courant's bracket to the space
$\Gamma(\mathcal{E}^1(M))$ of smooth sections of
$\mathcal{E}^1(M)$ and shows that the graph of the vector bundle
morphism $(\Lambda,E)^{\#} : T^{\ast}M\times \R \to TM\times \R$
is a Dirac subbundle of $\mathcal{E}^1(M)$. But the extended
bracket does not endow $\mathcal{E}^1(M)$ with a Courant algebroid
structure. A second approach of the problem is the one proposed by
the second author and J. Clemente-Gallardo in the recent paper
\cite{jj}. They introduce the notions of {\it generalized Courant
algebroid} (which is equivalent to the notion of {\it
Courant-Jacobi algebroid} independently defined by J. Grabowski
and G. Marmo in \cite{gm2}) and of {\it Dirac structure for a
generalized Courant algebroid} and give several connections
between Dirac structures for generalized Courant algebroids and
Jacobi manifolds. We note that the construction of \cite{jj}
includes as particular case the one of Wade and that the main
example of generalized Courant algebroid over $M$ is the direct
sum of a generalized Lie bialgebroid over $M$.

\vspace{1mm}

In the present work, by using the results mentioned above, we
establish a reduction theorem of Jacobi manifolds (Theorem
\ref{th-reduction}). It is well known that there are already
several geometric and algebraic treatments of the Jacobi reduction
problem (see, for instance, \cite{j2}, \cite{j3}, \cite{mk},
\cite{ib}). But, it is an original goal of the Dirac structures
theory to describe Jacobi reduction and to construct a more
general framework for the study of the related problems concerning
the projection of Jacobi structures and the existence of Jacobi
structures on certain submanifolds of Jacobi manifolds. Precisely,
on the way to our principal result, we construct an one to one
correspondence between Dirac subbundles, satisfying a certain
regularity condition, of the double $(TM\times \R)\oplus
(T^{\ast}M\times \R)$, where $M$ is a Jacobi manifold, and
quotient Jacobi manifolds of $M$ (Theorem \ref{basic-th}). Also,
the Reduction Theorem \ref{th-reduction} allows us to state
sufficient conditions under which a submanifold $N$ of
$(M,\Lambda,E)$ inherits a Jacobi structure, that include as
particular cases the results presented in \cite{i1}, \cite{dlm}.

\vspace{1mm}

The paper is organized as follows. In sections 2 and 3 we recall
some basic definitions and results concerning, respectively,
Jacobi structures, generalized Lie bialgebroids and Dirac
structures for generalized Courant algebroids. In section 4 we
establish a correspondence between Dirac structures and quotient
Jacobi manifolds (Theorem \ref{basic-th}). Using this
correspondence and the results for the pull-backs Dirac structures
by Lie algebroid morphisms, we prove, in section 5, a Jacobi
reduction theorem (Theorem \ref{th-reduction}) which is
essentially the Reduction Theorem proved in \cite{j2} and
independently in \cite{mk}. Finally, in section 6 we present some
applications and examples.

\vspace{3mm} \noindent {\bf Notation :} In this paper, $M$ is a
$C^{\infty}$-differentiable manifold of finite dimension. We
denote by $TM$ and $T^{\ast}M$, respectively, the tangent and
cotangent bundles over $M$, $C^{\infty}(M,\R)$ the space of all
real $C^{\infty}$-differentiable functions on $M$, $\Omega^k(M)$
the space of all differentiable $k$-forms on $M$ and
$\mathcal{X}(M)$ the space of all differentiable vector fields on
$M$. Also, we denote by $\delta$ the usual differential operator
on the graded space $\Omega(M)=\oplus_{k\in Z\!\!\!Z}\Omega^k(M)$.
For the Schouten bracket and the interior product of a form with a
multi-vector field, we use the convention of sign indicated by
Koszul \cite{kz}, (see also \cite{mrl}).

\section{Jacobi structures and Generalized Lie bialgebroids}
A {\it Jacobi manifold} is a differentiable manifold $M$ equipped
with a bivector field $\Lambda$ and a vector field $E$ such that
$$
[\Lambda,\Lambda]=-2E\wedge \Lambda
\hspace{5mm}\mathrm{and}\hspace{5mm}[E,\Lambda] = 0,
$$
where $[\,,\,]$ denotes the Schouten bracket (\cite{lch}). In this
case, $(\Lambda,E)$ defines on $C^{\infty}(M,\R)$ the internal
composition law $\{\,,\,\}_{(\Lambda,E)} : C^{\infty}(M,\R)\times
C^{\infty}(M,\R) \to C^{\infty}(M,\R)$ given, for all $f,g\in
C^{\infty}(M,\R)$, by
\begin{equation}\label{crJ}
\{f,g\}_{(\Lambda,E)} = \Lambda (\delta f,\delta g) + \langle
f\delta g - g \delta f, E\rangle,
\end{equation}
which endows $C^{\infty}(M,\R)$ with a local Lie algebra structure
\cite{kr}, \cite{lch}, (or with a Jacobi algebra structure in the
terminology of J. Grabowski {\it et al.} \cite{gr}, \cite{gm2}).

Let $(M_1,\Lambda_1,E_1)$ and $(M_2,\Lambda_2,E_2)$ be two Jacobi
manifolds and $\Psi : M_1 \to M_2$ a differentiable map. If
$\Lambda_1$ and $E_1$ are projectable by $\Psi$ on $M_2$ and their
projections are, respectively, $\Lambda_2$ and $E_2$, i.e
$\Psi_{\ast}\Lambda_1=\Lambda_2$ and $\Psi_{\ast}E_1=E_2$, then
$\Psi : M_1\to M_2$ is said to be a {\it Jacobi morphism} or a
{\it a Jacobi map}. When $\Psi : M_1\to M_2$ is a diffeomorphism,
the Jacobi structures $(\Lambda_1,E_1)$ and $(\Lambda_2,E_2)$ are
said to be {\it equivalent}.

\vspace{3mm}

A {\it Lie algebroid} over a smooth manifold $M$ is a vector
bundle $A\rightarrow M$ with a Lie algebra structure $[\,,\,]$ on
the space $\Gamma(A)$ of the global cross sections of $A\to M$ and
a bundle map $a : A \to TM$, called the {\it anchor map}, such
that
\begin{enumerate}
\item the homomorphism $a : (\Gamma(A),[\,,\,])
\to(\mathcal{X}(M),[\,,\,])$, induced by the anchor map, is a Lie
algebra homomorphism ; \item for all $f\in C^{\infty}(M,\R)$ and
for all $X,Y \in \Gamma(A)$,
$$
[X,fY] = f[X,Y] + (a(X)f)Y.
$$
\end{enumerate}
We denote a Lie algebroid over $M$ by the triple $(A,[\,,\,],a)$.
For more details see, for example, \cite{mck}, \cite{acw} and
\cite{mrl}.

\vspace{2mm}

A trivial example of a Lie algebroid over a differentiable
manifold $M$ is the triple $(TM\times \R, [\,,\,], \pi)$; for all
$(X,f),(Y,g)\in \Gamma(TM\times \R)\cong \mathcal{X}(M)\times
C^{\infty}(M,\R)$,
\begin{equation}\label{34}
[(X,f),(Y,g)] = ([X,Y],\, X\cdot g - Y\cdot f),
\end{equation}
and $\pi : TM\times \R \to TM$ is the projection on the first
factor.

\vspace{3mm}

The {\it Lie algebroid of a Jacobi manifold} $(M,\Lambda,E)$ is
defined in \cite{krb} as follows. We consider the vector bundle
$T^{\ast}M \times \R$ over $M$ and the vector bundle morphism
$(\Lambda,E)^{\#} : T^{\ast}M \times \R \to TM\times \R$ given,
for any $(\alpha,f) \in \Gamma(T^{\ast}M \times \R)$, by
$$
(\Lambda,E)^{\#}((\alpha,f)) = (\Lambda^{\#}(\alpha) +
fE,\,-\langle \alpha,E\rangle ).
$$
On the space $\Gamma(T^{\ast}M \times \R) \cong \Omega^1(M)\times
C^{\infty}(M,\R)$ we define the Lie algebra bracket
$[\,,\,]_{(\Lambda,E)}$ given, for all $(\alpha,f), (\beta,g) \in
\Gamma(T^{\ast}M \times \R)$, by
\begin{equation}\label{1}
[(\alpha,f), (\beta,g)]_{(\Lambda,E)} : = (\gamma,h),
\end{equation}
where
$$
\gamma : = \mathcal{L}_{\Lambda^{\#}(\alpha)}\beta -
\mathcal{L}_{\Lambda^{\#}(\beta)}\alpha -
\delta(\Lambda(\alpha,\beta)) + f\mathcal{L}_E\beta -
g\mathcal{L}_E\alpha - i_E(\alpha\wedge \beta),
$$
$$
h : = -\Lambda(\alpha,\beta) + \Lambda(\alpha,\delta g) -
\Lambda(\beta,\delta f) + \langle f\delta g - g\delta f, E\rangle.
$$
Then the triple $(T^{\ast}M \times \R, [\,,\,]_{(\Lambda,E)}, \pi
\circ (\Lambda,E)^{\#})$ is a Lie algebroid over $M$.

\vspace{3mm}

For a Lie algebroid $(A, [\,,\,], a)$ over $M$, we denote by
$A^{\ast}$ its dual vector bundle over $M$ and by $\bigwedge
A^{\ast} = \oplus_{k \in \Z}\bigwedge^k A^{\ast}$ the graded
exterior algebra of $A^{\ast}$. Sections of $\bigwedge A^{\ast}$
are called {\it A-differential forms} (or {\it A-forms}) on $M$.
There exists a graded endomorphism $d : \Gamma(\bigwedge A^{\ast})
\to \Gamma(\bigwedge A^{\ast})$, of degree 1, of the exterior
algebra of $A$-forms, called the {\it exterior derivative}, taking
an $A$-$k$-form $\eta$ to an $A$-$(k+1)$-form $d\eta$ such that,
for all $X_1, \ldots, X_{k+1} \in \Gamma(A)$,
\begin{eqnarray*}
d\eta(X_1, \ldots, X_{k+1}) & = &
\Sigma_{i=1}^{k+1}(-1)^{i+1}a(X_i)\cdot
\eta(X_1, \ldots, \hat{X}_i,\ldots, X_{k+1}) \\
& + & \Sigma_{1\leq i<j\leq k+1}(-1)^{i+j}\eta([X_i,X_j], X_1,
\ldots, \hat{X}_i, \ldots, \hat{X}_j, \ldots, X_{k+1}).
\end{eqnarray*}
The Lie algebroid axioms of $A$ imply that $d$ is a
$C^{\infty}(M,\R)$-multilinear superderivation of degree 1 such
that $d^2 = 0$. Also, we denote by $\bigwedge A = \oplus_{k \in
\Z}\bigwedge^k A$ the graded exterior algebra of $A$ whose
sections are called {\it A-multivector fields}. The Lie bracket on
$\Gamma(A)$ can be extended to the exterior algebra of
$A$-multivector fields and the result is a graded Lie bracket
$[\,,\,]$, called the {\it Schouten bracket} of the Lie algebroid
$A$. Details may be found, for instance, in \cite{mck}, \cite{ks1}
and \cite{acw}.

\vspace{3mm}

Let $(A, [\,,\,], a)$ be a Lie algebroid over $M$ and $\phi \in
\Gamma(A^{\ast})$ be an 1-cocycle in the Lie algebroid cohomology
complex with trivial coefficients (\cite{mck}, \cite{im1}), i.e.
for any $X,Y \in \Gamma(A)$,
\begin{equation}\label{2}
\langle \phi, [X,Y]\rangle = a(X)(\langle \phi, Y\rangle) -
a(Y)(\langle \phi, X\rangle).
\end{equation}
We modify the usual representation of the Lie algebra $(\Gamma(A),
[\,,\,])$ on the space $C^{\infty}(M,\R)$ by defining a new
representation $a^{\phi} : \Gamma(A)\times C^{\infty}(M,\R)\to
C^{\infty}(M,\R)$ as
\begin{equation}\label{3}
a^{\phi}(X,f) = a(X)f + \langle \phi,X\rangle f, \hspace{5mm}
\forall \, (X,f) \in \Gamma(A)\times C^{\infty}(M,\R).
\end{equation}
The resulting cohomology operator $d^{\phi} : \Gamma(\bigwedge
A^{\ast}) \to \Gamma(\bigwedge A^{\ast})$ of the new cohomology
complex is called the $\phi$-{\it differential} of $A$ and its
expression in terms of $d$ is
\begin{equation}\label{4}
d^{\phi}\eta = d\eta + \phi\wedge \eta, \hspace{5mm} \forall \,
\eta \in \Gamma(\bigwedge {}^k A^{\ast}).
\end{equation}
$d^{\phi}$ allows us to define, in a natural way, the $\phi$-{\it
Lie derivative by} $X\in \Gamma(A)$, $\mathcal{L}_X^{\phi} :
\Gamma(\bigwedge^k A^{\ast}) \to \Gamma(\bigwedge^k A^{\ast})$, as
the commutator of $d^{\phi}$ and of the contraction by $X$, i.e.
$\mathcal{L}_X^{\phi} = d^{\phi}\circ i_X + i_X \circ d^{\phi}$.
Its expression in terms of the usual Lie derivative $\mathcal{L}_X
= d\circ i_X + i_X\circ d$ is, for any $\eta \in
\Gamma(\bigwedge^k A^{\ast})$,
\begin{equation}\label{5}
\mathcal{L}_X^{\phi}\eta = \mathcal{L}_X\eta + \langle
\phi,X\rangle \eta.
\end{equation}
Using $\phi$ we can also modify the Schouten bracket $[\,,\,]$ on
$\Gamma(\bigwedge A)$ to the $\phi$-{\it Schouten bracket}
$[\,,\,]^{\phi}$ on $\Gamma(\bigwedge A)$. It is defined, for all
$P\in \Gamma(\bigwedge^p A)$ and $Q \in \Gamma(\bigwedge^q A)$, by
\begin{equation}\label{6}
[P,Q]^{\phi} = [P,Q] + (p-1)P\wedge (i_{\phi}Q) +
(-1)^p(q-1)(i_{\phi}P)\wedge Q,
\end{equation}
where $i_{\phi}Q$ can be interpreted as the usual contraction of a
multivector field with an 1-form. We remark that, when $p=q=1$,
$[P,Q]^{\phi} = [P,Q]$, i.e. the brackets $[\,,\,]^{\phi}$ and
$[\,,\,]$ coincide on $\Gamma(A)$. For a representation of the
differential calculus using the $\phi$-modified derivative, Lie
derivative and Schouten bracket, see \cite{im1} and \cite{gm1}.

\vspace{3mm}

The notion of {\it generalized Lie bialgebroid} has been
introduced by D. Iglesias and J.C. Marrero in \cite{im1} in such a
way that a Jacobi manifold has a generalized Lie bialgebroid
canonically associated and conversely. We consider a Lie algebroid
$(A, [\,,\,], a)$ over $M$ and an 1-cocycle $\phi \in
\Gamma(A^{\ast})$ and we assume that the dual vector bundle
$A^{\ast}\to M$ admits a Lie algebroid structure $([\,,\,]_{\ast},
a_{\ast})$ and that $W \in \Gamma(A)$ is an 1-cocycle in the Lie
algebroid cohomology complex with trivial coefficients of
$(A^{\ast},[\,,\,]_{\ast}, a_{\ast})$. Then, we say that :
\begin{definition}\label{D1}
The pair $((A,\phi),(A^{\ast},W))$ is a generalized Lie
bialgebroid over $M$ if, for all $X,Y\in \Gamma(A)$ and $P\in
\Gamma(\bigwedge^pA)$, the following conditions hold :
\begin{equation}\label{a1}
d_{\ast}^W[X,Y] = [d_{\ast}^WX,Y]^{\phi} + [X,d_{\ast}^WY]^{\phi}
\hspace{3mm} and \hspace{3mm} \mathcal{L}_{\ast \phi}^WP +
\mathcal{L}_W^{\phi}P = 0,
\end{equation}
where $d_{\ast}^W$ and $\mathcal{L}_{\ast}^W$ are, respectively,
the $W$-differential and the $W$-Lie derivative of $A^{\ast}$.
\end {definition}

An equivalent definition of this notion was presented in
\cite{gm1} by J. Grabowski and G. Marmo under the name of {\it
Jacobi bialgebroid}. Precisely, they define that :

\begin{definition}\label{D3}
The pair $((A,\phi),(A^{\ast},W))$ is a Jacobi bialgebroid if for
all $P\in \Gamma(\bigwedge^pA)$ and $Q\in \Gamma(\bigwedge^qA)$,
$$
d_{\ast}^W[P,Q]^{\phi} = [d_{\ast}^WP,Q]^{\phi} +
(-1)^{p+1}[P,d_{\ast}^WQ]^{\phi}.
$$
\end{definition}

In the particular case where $\phi = 0$ and $W=0$, by the above
two definitions we recover, respectively, the notion of {\it Lie
bialgebroid} introduced by K. Mackenzie and P. Xu in \cite{mx} and
its equivalent definition given by Yv. Kosmann-Schwarzbach in
\cite{ks1}.

\begin{remark}\label{R1}
{\rm The property of duality of a Lie bialgebroid is also verified
in the case of a generalized Lie bialgebroid : i.e. if
$((A,\phi),(A^{\ast},W))$ is a generalized Lie bialgebroid, so is
$((A^{\ast},W),(A,\phi))$ (see \cite{im1}, \cite{gm1}).
Consequently, conditions of Definition \ref{D1} as well as of
Definition \ref{D3} can be replaced by their dual versions.}
\end{remark}

The fundamental results of \cite{im1}, which will be used in the
sequel, are the following theorems.

\begin{theorem}\label{Jac-bialg}
Let $(M,\Lambda,E)$ be a Jacobi manifold. Then $\big((TM\times
\R,[\,,\,],\pi, (0,1)),(T^{\ast}M\times \R, [\,,\,]_{(\Lambda,E)},
\pi\circ (\Lambda,E)^{\#}, (-E,0))\big)$ is a generalized Lie
bialgebroid.
\end{theorem}

\begin{theorem}\label{bialg-Jac}
Let $((A,\phi),(A^{\ast},W))$ be a generalized Lie bialgebroid
over $M$. Then the bracket $\{\,,\,\}_ J : C^{\infty}(M,\R)\times
C^{\infty}(M,\R) \to C^{\infty}(M,\R)$ given, for all $f,g \in
C^{\infty}(M,\R)$, by
\begin{equation}\label{10}
\{f,g\}_J : = \langle d^{\phi}f, d_{\ast}^Wg\rangle,
\end{equation}
defines a Jacobi structure on $M$.
\end{theorem}

\begin{corol}\label{corol-bialg-Jac}
Let $\big((TM\times \R,[\,,\,],\pi, (0,1)),(T^{\ast}M\times \R,
[\,,\,]_{(\Lambda,E)}, \pi\circ (\Lambda,E)^{\#}, (-E,0))\big)$ be
the generalized Lie bialgebroid associated to a Jacobi manifold
$(M,\Lambda,E)$. Then,
\begin{equation}\label{7}
\hspace{35mm} \{f,g\}_J = \{f,g\}_{(\Lambda,E)}, \hspace{10mm}
\forall \,f,g \in C^{\infty}(M,\R).
\end{equation}
\end{corol}
{\bf Proof :} Effectively, for all $f,g \in C^{\infty}(M,\R)$,
\begin{eqnarray*}
\{f,g\}_J & \stackrel{(\ref{10})}{=} & \langle d^{(0,1)}f,
d_{\ast}^{(-E,0)}g\rangle = \langle (\delta f,f),
(-\Lambda^{\#}(\delta g)-gE, \langle \delta g,E\rangle)\rangle \\
          & = & \Lambda (\delta f, \delta g) + \langle f \delta g - g \delta f,
           E\rangle \stackrel{(\ref{crJ})}{=}  \{f,g\}_{(\Lambda,E)}. \,\bullet
\end{eqnarray*}

\vspace{3mm}

An important class of generalized Lie bialgebroids is the one of
{\it triangular generalized Lie bialgebroids} defined, also in
\cite{im1} and \cite{im3}, as follows :
\begin{definition}\label{trbi}
A generalized Lie bialgebroid $((A,\phi),(A^{\ast},W))$ is said to
be a triangular generalized Lie bialgebroid if there exists $P \in
\Gamma(\bigwedge^2A)$ such that $[P,P]^{\phi} = 0$, $a_{\ast} = a
\circ P^{\#}$, $W = -P^{\#}(\phi)$ and the Lie bracket
$[\,,\,]_{\ast}$ on $\Gamma(A^{\ast})$ is the bracket
\begin{equation}\label{P-br}
[\alpha, \beta]_P = \mathcal{L}^{\phi}_{P^{\#}(\alpha)}\beta -
\mathcal{L}^{\phi}_{P^{\#}(\beta)}\alpha
-d^{\phi}(P(\alpha,\beta)), \hspace{10mm} \forall \, \alpha,\beta
\in \Gamma(A^{\ast}).
\end{equation}
\end{definition}

A characteristic example of triangular generalized Lie bialgebroid
is the generalized Lie bialgebroid of a Jacobi manifold
$(M,\Lambda,E)$ (Theorem \ref{Jac-bialg}), where
$[(\Lambda,E),(\Lambda,E)]^{(0,1)}=0$ holds.

\section{Generalized Courant algebroids and Dirac structures}

The notion of {\it generalized Courant algebroid} has been
introduced by the second author and J. Clemente-Gallardo in
\cite{jj} and independently, under the name of {\it Courant-Jacobi
algebroid}, by J. Grabowski and G. Marmo in \cite{gm2}. In this
section, we recall some basic facts concerning this notion and its
relation with Dirac and Jacobi structures.

\begin{definition}[\cite{jj}]\label{Def-GenCour}
A generalized Courant algebroid over a smooth manifold $M$ is a
vector bundle $E$ over $M$ equipped with : (i) a nondegenerate
symmetric bilinear form $(\,,\,)$ on the bundle, (ii) a
skew-symmetric bracket $[\,,\,]$ on $\Gamma(E)$, (iii) a bundle
map $\rho : E \to TM$ and (iv) an $E$-1-form $\theta$ such that,
for any $e_1,e_2 \in \Gamma(E)$, $\langle \theta, [e_1,e_2]\rangle
= \rho(e_1)\langle \theta,e_2\rangle - \rho(e_2)\langle
\theta,e_1\rangle$, which satisfy the following relations :
\begin{enumerate}
\item for any $e_1,e_2,e_3 \in \Gamma(E)$,
$$
[[e_1,e_2],e_3] + c.p. = \mathcal{D}^{\theta}T(e_1,e_2,e_3) ;
$$
\item for any $e_1,e_2 \in \Gamma(E)$,
\begin{equation}\label{2def}
\rho([e_1,e_2]) = [\rho(e_1), \rho(e_2)] ;
\end{equation}
\item for any $e_1,e_2 \in \Gamma(E)$ and $f\in C^{\infty}(M,\R)$,
\begin{equation}\label{3def}
[e_1,fe_2] = f[e_1, e_2] + (\rho(e_1)f)e_2 -
(e_1,e_2)\mathcal{D}f;
\end{equation}
\item for any $f,g \in C^{\infty}(M,\R)$,
$$
(\mathcal{D}^{\theta}f,\mathcal{D}^{\theta}g) = 0 ;
$$
\item for any $e,e_1,e_2 \in \Gamma(E)$,
$$ \rho(e)(e_1,e_2) +
\langle \theta, e\rangle(e_1,e_2) = ([e,e_1]
+\mathcal{D}^{\theta}(e,e_1),e_2) + (e_1,[e,e_2]
+\mathcal{D}^{\theta}(e,e_2)).
$$
\end{enumerate}
For any $e_1,e_2,e_3 \in \Gamma(E)$, $T(e_1,e_2,e_3)$ is the
function on the base $M$ defined by
$$
T(e_1,e_2,e_3) = \frac{1}{3}([e_1,e_2],e_3) + c.p.
$$
$\mathcal{D}, \mathcal{D}^{\theta} : C^{\infty}(M,\R) \to
\Gamma(E)$ are the maps defined, for any $f \in C^{\infty}(M,\R)$,
by $\mathcal{D}f = \frac{1}{2}\beta^{-1}\rho^{\ast}\delta f$ and
$\mathcal{D}^{\theta}f = \mathcal{D}f +
\frac{1}{2}f\beta^{-1}(\theta)$, $\beta$ being the isomorphism
between $E$ and $E^{\ast}$ defined by the nondegenerate bilinear
form $(\,,\,)$. In other words, for any $e\in \Gamma(E)$,
$$
(\mathcal{D}f,e) = \frac{1}{2}\rho(e)f \hspace{4mm}and\hspace{4mm}
(\mathcal{D}^{\theta}f,e) = \frac{1}{2}(\rho(e)f + \langle
\theta,e\rangle f).
$$
\end{definition}

The above definition is based on the original definition of
Courant algebroid presented in \cite{lwx1} by Z.-J. Liu {\it et
al.} while its equivalent definition proposed in \cite{gm2} is
based on the alternative definition of Courant algebroid given by
D. Roytenberg in \cite{ro}. Their equivalence is established in
\cite{jj}.

\vspace{1mm} By defining, for any $e\in \Gamma(E)$, the first
order differential operator $\rho^{\theta}(e)$ by
\begin{equation}\label{def-rho-theta}
\rho^{\theta}(e)=\rho(e)+\langle \theta,e \rangle,
\end{equation}
we have that (\ref{2def}) is equivalent (\cite{jj}) to
\begin{equation}\label{33}
\rho^{\theta}([e_1,e_2])=[\rho^{\theta}(e_1),\rho^{\theta}(e_2)],
\end{equation}
where the bracket on the right-hand side is the Lie bracket
(\ref{34}).

\begin{definition}\label{D2}
A Dirac structure for a generalized Courant algebroid $(E,\theta)$
over $M$ is a subbundle $L\subset E$ that is maximal isotropic
under $(\,,\,)$ and integrable, i.e. $\Gamma(L)$ is closed under
$[\,,\,]$.
\end{definition}

It is immediate from the above definition that a Dirac subbundle
$L$ of $(E,\theta)$ is a Lie algebroid under the restrictions of
the bracket $[\,,\,]$ and of the anchor $\rho$ to $\Gamma(L)$. If
$\theta \in \Gamma(L^{\ast})$, then it is an 1-cocycle for the Lie
algebroid cohomology with trivial coefficients of
$(L,[\,,\,]\vert_{L},\rho\vert_L)$.

\vspace{3mm}

We consider now a generalized Lie bialgebroid
$((A,\phi),(A^{\ast},W))$ over $M$ and we denote by $E$ its vector
bundle direct sum, i.e. $E = A\oplus A^{\ast}$. On $E$ there exist
two natural nondegenerate bilinear forms, one symmetric
$(\,,\,)_+$ and another skew-symmetric $(\,,\,)_-$ given, for any
$e_1 = X_1+\alpha_1, e_2 = X_2+\alpha_2 \in E = A\oplus A^{\ast}$,
by
\begin{equation}\label{13}
(e_1,e_2)_{\pm} = (X_1+\alpha_1,X_2+\alpha_2)_{\pm} =
\frac{1}{2}(\langle \alpha_1,X_2\rangle \pm \langle
\alpha_2,X_1\rangle).
\end{equation}
On $\Gamma(E)$, which is identified with $\Gamma(A)\oplus
\Gamma(A^{\ast})$, we introduce the bracket $\lcf \,,\,\rcf$
defined, for all $e_1 = X_1 + \alpha_1, e_2 = X_2 + \alpha_2 \in
\Gamma(E)$, by
\begin{eqnarray}\label{14}
\lcf e_1,e_2 \rcf & \!\!\! = \!\!\! & \lcf X_1 + \alpha_1, X_2 + \alpha_2 \rcf   \nonumber\\
                  & \!\!\! = \!\!\! & ([X_1,X_2]^{\phi} + \mathcal{L}_{\ast \alpha_1}^WX_2 -
                         \mathcal{L}_{\ast \alpha_2}^WX_1 - d_{\ast}^W(e_1,e_2)_- ) + \nonumber \\
                  & \!\!\!  \!\!\! & ([\alpha_1,\alpha_2]_{\ast}^W + \mathcal{L}_{X_1}^{\phi}\alpha_2 -
                         \mathcal{L}_{X_2}^{\phi}\alpha_1 +
                         d^{\phi}(e_1,e_2)_-).
\end{eqnarray}
Finally, let $\rho : E \to TM$ be the bundle map given by $\rho =
a + a_{\ast}$, i.e., for any $X+\alpha \in E$,
\begin{equation}\label{15}
\rho(X+\alpha) = a(X) + a_{\ast}(\alpha).
\end{equation}
The following result, which is proved in \cite{jj}, shows that the
notion of generalized Courant algebroid permits us to generalize
the double construction for Lie bialgebras (the {\it Drinfeld
double}, \cite{dr}) and Lie bialgebroids (\cite{lwx1}) to
generalized Lie bialgebroids.

\begin{theorem}[\cite{jj}]
If $((A,\phi),(A^{\ast},W))$ is a generalized Lie bialgebroid over
$M$, then $E = A \oplus A^{\ast}$ endowed with $(\lcf \,,\,\rcf,
(\,,\,)_+,\rho)$ and $\theta = \phi + W \in \Gamma(E^{\ast})$ is a
generalized Courant algebroid over $M$. The operators
$\mathcal{D}$ and $\mathcal{D}^{\theta}$ are, respectively,
$\mathcal{D} = (d_{\ast}+d)\vert_{C^{\infty}(M,I \!\! R)}$ and
$\mathcal{D}^{\theta} = (d_{\ast}^W +
d^{\phi})\vert_{C^{\infty}(M,I\!\!R)}$.
\end{theorem}

There are two important classes of Dirac structures for the
generalized Courant algebroid $(E,\theta) = (A\oplus A^{\ast},
\phi +W)$ studied in \cite{jj}.

\vspace{3mm} \noindent {\bf The Dirac structure of the graph of an
{\it A}-bivector field :} Let $\Omega$ be an $A$-bivector field
and $\Omega^{\#} : A^{\ast} \to A$ the associated vector bundle
map. The graph of $\Omega^{\#}$ is the maximal isotropic vector
subbundle
$$
L = \{\Omega^{\#}\alpha + \alpha \, / \, \alpha \in A^{\ast}\}
$$
of $(A\oplus A^{\ast},(\,,\,)_+)$. $L$ is a Dirac structure for
$(A\oplus A^{\ast}, \phi +W)$ if and only if $\Omega$ satisfies
the Maurer-Cartan type equation :
$$
d_{\ast}^W\Omega + \frac{1}{2}[\Omega,\Omega]^{\phi} = 0.
$$

\vspace{3mm} \noindent {\bf Null Dirac structures :} Let $D
\subset A$ be a vector subbundle of $A$ and $D^{\bot}\subset
A^{\ast}$ its conormal bundle, i.e.
\begin{equation}\label{16}
D^{\bot} = \{\alpha \in A^{\ast}\,/\, \langle \alpha,X\rangle =0,
\hspace{2mm} \forall X \in A \}.
\end{equation}
Then, $L = D\oplus D^{\bot}$ is a Dirac structure for $(A\oplus
A^{\ast}, \phi +W)$ if and only if $D$ and $D^{\bot}$ are Lie
subalgebroids (\cite{mck}) of $A$ and $A^{\ast}$, respectively. It
is clear that in this context, as in the context of a Lie
bialgebroid, $L = D\oplus D^{\bot}$ if and only if the
skew-symmetric bilinear form $(\,,\,)_-$, defined on $E=A\oplus
A^{\ast}$ by (\ref{13}), vanishes on $L$. For this reason, $L$ is
said to be a {\it null Dirac structure}.

\vspace{3mm} A third important category of Dirac structures for
$(E,\theta) = (A\oplus A^{\ast}, \phi +W)$, also studied in
\cite{jj}, which generalizes both the above presented categories,
is :

\vspace{3mm} \noindent {\bf Dirac structures defined by a
characteristic pair :} We consider a pair $(D,\Omega)$ of a smooth
subbundle $D \subset A$ and of an $A$-bivector field $\Omega$. We
construct, following \cite{liu}, a subbunlde $L\subset A\oplus
A^{\ast}$ by setting :
\begin{equation}\label{17}
L = \{X + \Omega^{\#}\alpha + \alpha\,/\, X \in D \hspace{2mm}and
\hspace{2mm} \alpha \in D^{\bot}\} = D\oplus graph
(\Omega^{\#}\vert_{D^{\bot}}).
\end{equation}
$L$ is maximal isotropic with respect to $(\,,\,)_+$. The pair
$(D,\Omega)$ is called the {\it characteristic pair of L} while
the subbundle $D = L\cap (A\oplus \{0\})$, also denoted by $D = L
\cap A$, is called the {\it characteristic subbundle of L}.

\vspace{1mm} For simplicity, we will assume in the sequel that $D
= L\cap A$ is of constant rank.

\vspace{1mm}

Moreover, since $D^{\bot}$ may be considered as the dual bundle
$(A/D)^{\ast}$ of the quotient bundle $A/D$, the restricted vector
bundle map $\Omega^{\#}\vert_{D^{\bot}}$ can be seen as the bundle
map associated to an $A/D$-bivector field. Hence, two pairs
$(D_1,\Omega_1)$ and $(D_2,\Omega_2)$ of a smooth subbundle and of
an $A$-bivector field determine the same subbundle $L\subset
A\oplus A^{\ast}$ (given by (\ref{17})) if and only if
\begin{equation}\label{18}
D_1 = D_2 = : D \hspace{3mm}\mathrm{and}\hspace{3mm}
\Omega_1^{\#}(\alpha) - \Omega_2^{\#}(\alpha)\in D, \hspace{3mm}
\forall \alpha \in D^{\bot}.
\end{equation}
Let $pr : \Gamma(\bigwedge A)\to \Gamma(\bigwedge (A/D))$ be the
map on the spaces of sections, induced by the natural projection
$A\to A/D$. In order to express that the projection under $pr$ of
an $A$-multivector field $\Omega \in \Gamma(\bigwedge A)$ vanishes
in $\Gamma(\bigwedge (A/D))$, we write $\Omega \equiv 0(mod\,D)$.
Thus, the second condition of (\ref{18}) can be written as
$\Omega_1-\Omega_2 \equiv 0 (mod\,D)$.

\vspace{1mm}

The conditions under which $L=D\oplus graph
(\Omega^{\#}\vert_{D^{\bot}})$ is a Dirac subbundle of $(A\oplus
A^{\ast}, \phi +W)$ are given by :

\begin{theorem}[\cite{jj}]\label{charpair}
Let $L=D\oplus graph (\Omega^{\#}\vert_{D^{\bot}})$ be a maximal
isotropic subbundle of $A\oplus A^{\ast}$. Then, $L$ is a Dirac
structure for the generalized Courant algebroid $(A\oplus
A^{\ast}, \phi + W)$ if and only if
\begin{enumerate}
\item[i)] $D$ is a Lie subalgebroid of $A$ ; \item[ii)]
$d_{\ast}^W\Omega + \frac{1}{2}[\Omega,\Omega]^{\phi}\equiv
0(mod\,D)$ ; \item[iii)] $D^{\bot}$ is integrable for the sum
bracket $[\,,\,]_{\ast} + [\,,\,]_{\Omega}$, i.e., for all
$\alpha,\beta \in \Gamma(D^{\bot})$, $[\alpha,\beta]_{\ast} +
[\alpha,\beta]_{\Omega}\in \Gamma(D^{\bot})$, where
$[\,,\,]_{\Omega}$ is the bracket determined on $\Gamma(A^{\ast})$
by (\ref{P-br}).
\end{enumerate}
\end{theorem}

In the particular case where $((A,\phi),(A^{\ast},W))$ is a
triangular generalized Lie bialgebroid (Definition \ref{trbi}),
Theorem \ref{charpair} takes the following form :

\begin{corol}[\cite{jj}]\label{corol-charpair}
Let $((A,\phi),(A^{\ast},W),P)$ be a triangular generalized Lie
bialgebroid and $L\subset A\oplus A^{\ast}$, $L=D\oplus graph
(\Omega^{\#}\vert_{D^{\bot}})$, a maximal isotropic subbundle of
$A\oplus A^{\ast}$ with a fixed characteristic pair $(D,\Omega)$.
Then $L$ is a Dirac structure for the generalized Courant
algebroid $(A\oplus A^{\ast}, \phi + W)$ if and only if
\begin{enumerate}
\item[i)] $D$ is a Lie subalgebroid of $A$ ; \item[ii)] $[P+
\Omega,P+\Omega]^{\phi}\equiv 0(mod\,D)$ ; \item[iii)] for any $X
\in \Gamma(D)$, $\mathcal{L}_X^{\phi}(P+\Omega)\equiv 0(mod\,D)$.
\end{enumerate}
\end{corol}

\section{Jacobi structures and Dirac reducible subbundles}
We consider a generalized Lie bialgebroid
$((A,[\,,\,],a,\phi),(A^{\ast},[\,,\,]_{\ast},a_{\ast},W))$ over
$M$ and we construct the associated generalized Courant algebroid
$(A\oplus A^{\ast}, \lcf \,,\,\rcf,(\,,\,)_+, \rho, \theta)$ over
$M$, i.e. $\lcf\,,\,\rcf$ is determined by (\ref{14}), $\rho = a +
a_{\ast}$ and $\theta = \phi + W$. We introduce the notions of
{\it reducible Dirac structure} for $(A\oplus A^{\ast}, \lcf
\,,\,\rcf,(\,,\,)_+,\rho, \theta)$ and of {\it admissible
function} of a Dirac structure for $(A\oplus A^{\ast}, \lcf
\,,\,\rcf,(\,,\,)_+,\rho, \theta)$ in an analog manner as in the
case of a Dirac structure for a Lie bialgebroid (\cite{lwx2}).

\begin{definition}\label{Def-red}
We say that a Dirac subbundle $L$ for $(A\oplus A^{\ast}, \lcf
\,,\,\rcf,(\,,\,)_+,\rho, \theta)$ is reducible if the image
$a(D)$ of its characteristic subbundle $D=L\cap A$ by the anchor
map $a$ defines a simple foliation $\mathcal{F}$ of $M$. By the
term "simple foliation", we mean that $\mathcal{F}$ is a regular
foliation such that the space $M/\mathcal{F}$ is a nice manifold
and the canonical projection $M \to M/\mathcal{F}$ is a
submersion.
\end{definition}

\begin{definition}\label{Def-admis}
Let $L$ be a Dirac subbundle for $(A\oplus A^{\ast}, \lcf
\,,\,\rcf,(\,,\,)_+,\rho, \theta)$. We say that a function $f \in
C^{\infty}(M,\R)$ is $L$-admissible if there exists $Y_f \in
\Gamma(A)$ such that $Y_f + d^{\phi}f \in \Gamma(L)$.
\end{definition}
Obviously, $Y_f$ is unique up to a smooth section of $L\cap A$. We
denote by $C_L^{\infty}(M,\R)$ the set of all $L$-admissible
functions of $C^{\infty}(M,\R)$.

\vspace{3mm}

Let $L\subset A\oplus A^{\ast}$ be a Dirac structure for $(A\oplus
A^{\ast}, \lcf \,,\,\rcf,(\,,\,)_+, \rho, \theta)$. On
$C_L^{\infty}(M,\R)$ we define the bracket $\{\,,\,\}_L$ by
setting, for all $f,g \in C_L^{\infty}(M,\R)$,
\begin{equation}\label{19}
\{f,g\}_L : = \rho^{\theta}(e_f)g,
\end{equation}
where $e_f = Y_f + d^{\phi}f$. An equivalent expression of
(\ref{19}) is :
\begin{equation}\label{20}
\{f,g\}_L = \langle Y_f,d^{\phi}g\rangle + \{f,g\}_J,
\end{equation}
where $\{\,,\,\}_J$ is the bracket (\ref{10}) of the Jacobi
structure on $M$ defined by the generalized Lie bialgebroid
structure $((A,\phi),(A^{\ast},W))$ over $M$. Effectively,
\begin{eqnarray*}
\{f,g\}_L & \stackrel{(\ref{19})}{=} & \rho^{\theta}(e_f)g
          = ((a^{\phi} + a_{\ast}^W)(Y_f + d^{\phi}f))g
        = a^{\phi}(Y_f)g + a_{\ast}^W(d^{\phi}f)g \\
        &=& \langle Y_f,d^{\phi}g\rangle + \langle d^{\phi}f,
              d_{\ast}^Wg\rangle
         \stackrel{(\ref{10})}{=}  \langle Y_f,d^{\phi}g\rangle + \{f,g\}_J.
\end{eqnarray*}
It is easy to check that (\ref{19}) or equivalently (\ref{20}) is
well defined. In fact, if $Y'_f=Y_f+X$, with $X\in \Gamma(L\cap
A)$, is an other section of $A$ such that $e'_f=Y'_f+d^{\phi}f \in
\Gamma(L)$, we have
$$
\langle Y'_f,d^{\phi}g\rangle + \{f,g\}_J = \langle Y_f
+X,d^{\phi}g\rangle + \{f,g\}_J = \{f,g\}_L + \langle
X,d^{\phi}g\rangle = \{f,g\}_L,
$$
since, $L$ being isotropic, $(X+0,Y_g+d^{\phi}g)_+ = 0
\Leftrightarrow \langle X,d^{\phi}g\rangle = 0$.

\begin{proposition}\label{pr-admis}
The space $C_L^{\infty}(M,\R)$ endowed with the bracket
$\{\,,\,\}_L$, given by (\ref{19}), is a Lie algebra.
\end{proposition}
{\bf Proof :} We must prove that $C_L^{\infty}(M,\R)$ is closed
under $\{\,,\,\}_L$ and that $\{\,,\,\}_L$ is a bilinear,
skew-symmetric bracket which satisfies the Jacobi identity.

\vspace{2mm} \noindent {\it Closeness of $\{\,,\,\}_L$ in
$C_L^{\infty}(M,\R)$} : Let $f,g \in C_L^{\infty}(M,\R)$ be two
$L$-admissible functions. Then, there exist $Y_f,Y_g \in
\Gamma(A)$ such that $e_f = Y_f + d^{\phi}f, e_g = Y_g + d^{\phi}g
\in \Gamma(L)$. We consider the bracket $\lcf e_f,e_g\rcf$ ;
according to (\ref{14}), its component in $\Gamma(A^{\ast})$ is :
$$
[d^{\phi}f,d^{\phi}g]_{\ast}^W + \mathcal{L}_{Y_f}^{\phi}d^{\phi}g
- \mathcal{L}_{Y_g}^{\phi}d^{\phi}f + d^{\phi}(e_f,e_g)_-.
$$
We have (\cite{im1}),
$$
[d^{\phi}f,d^{\phi}g]_{\ast}^W  =  -
\mathcal{L}_{d_{\ast}^Wf}^{\phi}d^{\phi}g = - d^{\phi}\langle
d_{\ast}^Wf,d^{\phi}g\rangle = - d^{\phi}\{g,f\}_J =
d^{\phi}\{f,g\}_J
$$
and, on the other hand,
$$
\mathcal{L}_{Y_f}^{\phi}d^{\phi}g -
\mathcal{L}_{Y_g}^{\phi}d^{\phi}f + d^{\phi}(e_f,e_g)_- = -
d^{\phi}(e_f,e_g)_- = - d^{\phi}(e_f,e_g)_- +
\underbrace{d^{\phi}(e_f,e_g)_+}_{= 0} = d^{\phi}\langle
Y_f,d^{\phi}g\rangle.
$$
Thus,
$$
[d^{\phi}f,d^{\phi}g]_{\ast}^W + \mathcal{L}_{Y_f}^{\phi}d^{\phi}g
- \mathcal{L}_{Y_g}^{\phi}d^{\phi}f + d^{\phi}(e_f,e_g)_- =
d^{\phi}\{f,g\}_J + d^{\phi}\langle Y_f,d^{\phi}g\rangle
\stackrel{(\ref{20})}{=}d^{\phi}\{f,g\}_L,
$$
which means that $\{f,g\}_L$ is an $L$-admissible function, i.e.
$\{f,g\}_L\in C_L^{\infty}(M,\R)$, and that we can take $\lcf
e_f,e_g\rcf = e_{\{f,g\}_L}$.

\vspace{2mm} \noindent {\it Bilinearity and skew-symmetry of
$\{\,,\,\}_L$} : It is obvious that $\{\,,\,\}_L$ is bilinear.
Also, for any $f\in C_L^{\infty}(M,\R)$, we have $(e_f,e_f)_+ = 0
\Leftrightarrow \langle Y_f, d^{\phi}f\rangle = 0$, so $\{f,f\}_L
\stackrel{(\ref{20})}{=} \langle Y_f, d^{\phi}f\rangle + \{f,f\}_J
= 0+0 = 0$, which implies the skew-symmetry of $\{\,,\,\}_L$.

\vspace{2mm} \noindent {\it Jacobi identity} : By a
straightforward, but long, calculation we get that, for any
$f,g,h\in C_L^{\infty}(M,\R)$, the Jacobi identity holds :
$$
\{f,\{g,h\}_L\}_L + \{g,\{h,f\}_L\}_L + \{h,\{f,g\}_L\}_L = 0.
$$
Hence, $(C_L^{\infty}(M,\R),\{\,,\,\}_L)$ is a Lie algebra.
$\bullet$

\vspace{3mm}

In the particular case where the constant function $1$ is an
$L$-admissible function, $C_L^{\infty}(M,\R)$ equipped with the
usual product of functions $"\cdot"$ is an associative commutative
algebra with unit and $\{\,,\,\}_L$ is a first order differential
operator on each of its arguments. In fact, $1\in
C_L^{\infty}(M,\R)$ means that there exists $Y_1\in \Gamma(A)$
such that $Y_1+d^{\phi}1 = Y_1+\phi \in \Gamma(L)$. Then, for any
$f,g\in C_L^{\infty}(M,\R)$, $f\cdot g \in C_L^{\infty}(M,\R)$
since, for $Y_{fg} = fY_g+gY_f-fgY_1 \in \Gamma(A)$,
$Y_{fg}+d^{\phi}(fg) \in \Gamma(L)$. Moreover, for any $f,g,h\in
C_L^{\infty}(M,\R)$,
\begin{eqnarray*}
\{f,gh\}_L & \stackrel{(\ref{20})}{=} & \langle Y_f,
           d^{\phi}(gh)\rangle + \{f,gh\}_J \\
         & = & \langle Y_f, gdh + hdg + gh\phi \rangle +
         g\{f,h\}_J + h\{f,g\}_J - gh\{f,1\}_J \\
         & = & g(\langle Y_f, d^{\phi}h \rangle + \{f,h\}_J) +
               h(\langle Y_f, d^{\phi}g \rangle + \{f,g\}_J) -
               gh(\langle Y_f, \phi \rangle + \{f,1\}_J) \\
         & = & g\{f,h\}_L + h\{f,g\}_L - gh\{f,1\}_L
\end{eqnarray*}
and by the skew-symmetry of $\{\,,\,\}_L$ we obtain the desired
result. Consequently,

\begin{theorem}\label{th-admis}
If 1 is an $L$-admissible function, then
$(C_L^{\infty}(M,\R),\{\,,\,\}_L)$ is a Jacobi algebra.
\end{theorem}

The above result generalizes the one of A. Wade (\cite{wd}) for
the $\mathcal{E}^1(M)$-Dirac structures.

\vspace{3mm}

In the following we will establish the characteristic equations of
a Dirac structure (see, also, \cite{c}).

\begin{lemma}\label{charact-eq}
Let $L$ be a Dirac structure for $(A\oplus A^{\ast}, \lcf
\,,\,\rcf,(\,,\,)_+, \rho, \theta)$, $\varpi : A\oplus A^{\ast}
\to A$ and $\varpi_{\ast} : A\oplus A^{\ast} \to A^{\ast}$ the
natural projections from $A\oplus A^{\ast}$ onto $A$ and
$A^{\ast}$, respectively. Then, $\ker \varpi\vert_L = L\cap
A^{\ast}$ and $\ker \varpi_{\ast}\vert_L = L\cap A$. Also,
\begin{equation}\label{eq5}
\varpi(L)^{\bot} = L\cap A^{\ast} \hspace{3mm} and \hspace{3mm}
\varpi_{\ast}(L) = (L\cap A)^{\bot}.
\end{equation}
\end{lemma}
{\bf Proof :} We denote by $L_x$, $A_x$ and $A_x^{\ast}$ the
fibers over $x\in M$ of $L$, $A$ and $A^{\ast}$, respectively. It
is clear that, at each point $x\in M$, $\ker \varpi\vert_{L_x} =
L_x\cap A_x^{\ast}$ and $\ker \varpi_{\ast}\vert_{L_x} = L_x\cap
A_x$, thus $\dim (\ker \varpi\vert_{L_x}) = \dim (L_x\cap
A_x^{\ast})$ and $\dim (\ker \varpi_{\ast}\vert_{L_x}) = \dim
(L_x\cap A_x)$. Also, at each point $x \in M$,
$\varpi(L_x)^{\bot}=L_x\cap A_x^{\ast}$ and $\varpi_{\ast}(L_x) =
(L_x\cap A_x)^{\bot}$. Effectively, if $\alpha(x) \in L_x\cap
A_x^{\ast}$, then $0+\alpha (x) \in L_x$, thus, for any
$Y(x)+\beta (x) \in L_x$, $(0+\alpha (x), Y(x)+\beta (x))_+ = 0
\Leftrightarrow \langle \alpha (x), Y(x)\rangle = 0$ that implies
that $\alpha (x) \in \varpi(L_x)^{\bot}$, i.e $L_x\cap A_x^{\ast}
\subseteq \varpi(L_x)^{\bot}$, and by a dimension count we
conclude the equality. Analogously, we prove the second equation
of (\ref{eq5}) for the fibers over $x$. Since the above results
hold at each point $x\in M$, we get that the characteristic
equations (\ref{eq5}) of $L$ hold. $\bullet$

\begin{lemma}\label{lemma-1-admis}
The constant function $1$ is an $L$-admissible function if and
only if, for any $Y\in \Gamma(D)$,
\begin{equation}\label{eq1}
\langle \phi, Y\rangle = 0.
\end{equation}
\end{lemma}
{\bf Proof :} In fact, if $1\in C_L^{\infty}(M,\R)$, then there
exists $Y_1 \in \Gamma(A)$ such that $Y_1 + d^{\phi}1=Y_1+\phi \in
\Gamma(L)$. Also, for every $Y\in \Gamma(D)$, $Y+0 \in \Gamma(L)$
and $(L, (\,,\,)_+)$ is maximal isotropic. Thus, for any $Y\in
\Gamma(D)$, $(Y+0,Y_1+\phi)_+ = 0 \Leftrightarrow \langle \phi,
Y\rangle = 0$. Conversely, we suppose that, for any $Y\in
\Gamma(D)$, $\langle \phi, Y\rangle = 0$; then we will prove that
$1 \in C_L^{\infty}(M,\R)$. Effectively, if $1$ is not an
$L$-admissible function, then, for any $Y_1\in \Gamma(A)$, $Y_1+
d^{\phi}1=Y_1+\phi$ is not a section of $L$, fact which implies
that $\phi$ is not a section of
$\varpi_{\ast}(L)\stackrel{(\ref{eq5})}{=}(L\cap A)^{\bot} =
D^{\bot}$. Therefore, there exists $Y\in \Gamma(D)$ such that
$\langle \phi,Y\rangle \neq 0$; contradiction. $\bullet$

\begin{proposition}\label{prop-admis}
Let $L\subset A\oplus A^{\ast}$ be a reducible Dirac structure for
$(A\oplus A^{\ast}, \lcf \,,\,\rcf,(\,,\,)_+, \rho, \theta)$ and
$\mathcal{F}$ the simple foliation of $M$ defined by the
distribution $a(D)$, $D=L\cap A$, on $M$. If $1$ is an
$L$-admissible function, then $f\in C_L^{\infty}(M,\R)$ if and
only if $f$ is constant along the leaves of $\mathcal{F}$.
\end{proposition}
{\bf Proof :} Let $f$ be an $L$-admissible function, i.e. there
exists $Y_f \in \Gamma(A)$ such that $Y_f+d^{\phi}f \in
\Gamma(L)$, and $X \in \Gamma(a(D))$ a section of the distribution
$a(D)$. $X\in \Gamma(a(D))$ means that there exists $Y \in
\Gamma(D)$ such that $X=a(Y)$ and $Y\in \Gamma(D)$ means that $Y+0
\in \Gamma(L)$. Since $L$ is maximally isotropic,
$$
(Y+0, Y_f+d^{\phi}f)_+ = 0 \Leftrightarrow \langle d^{\phi}f,Y
\rangle = 0 \Leftrightarrow \langle df, Y\rangle + f \langle \phi,
Y\rangle = 0 \stackrel{(\ref{eq1})}{\Leftrightarrow} \langle
\delta f, a(Y)\rangle = 0 \Leftrightarrow \langle \delta f, X
\rangle =0.
$$
By the last equation, which holds for any $X\in \Gamma(a(D))$, we
get that $f$ is constant along the leaves of $\mathcal{F}$.
Conversely, let $f$ be a function on $M$ constant along the leaves
of $\mathcal{F}$, i.e., for any $X\in \Gamma(a(D))$, $\langle
\delta f, X \rangle =0$. But, $X\in \Gamma(a(D))$ means that there
exists $Y\in \Gamma(D)$ such that $X=a(Y)$. Thus, for any $X\in
\Gamma(a(D))$, $X=a(Y)$ with $Y\in \Gamma(D)$,
\begin{equation}\label{eq2}
\langle \delta f, X \rangle =0 \Leftrightarrow \langle \delta f,
a(Y)\rangle =0 \Leftrightarrow \langle df, Y\rangle = 0
\stackrel{(\ref{eq1})}{\Leftrightarrow} \langle df +f \phi,
Y\rangle = 0 \Leftrightarrow \langle d^{\phi}f, Y\rangle = 0,
\end{equation}
for any $Y\in \Gamma(D)$. If $f$ is not an $L$-admissible
function, then, for any $Z\in \Gamma(A)$, $Z+d^{\phi}f$ is not a
section of $L$. So, $d^{\phi}f$ is not a section of
$\varpi_{\ast}(L)\stackrel{(\ref{eq5})}{=}(L\cap A)^{\bot} =
D^{\bot}$. Therefore, there exists $Y\in \Gamma(D)$ such that
$\langle d^{\phi}f,Y\rangle \neq 0$; contradiction. $\bullet$

\vspace{3mm}

By the above study we conclude :

\begin{theorem}\label{th-dir-jacquot}
Let $L$ be a reducible Dirac subbundle for $(A\oplus A^{\ast},
\lcf \,,\,\rcf,(\,,\,)_+, \rho, \theta)$. We suppose that $1$ is
an $L$-admissible function. Then $L$ induces a Jacobi structure on
$M/\mathcal{F}$ defined by the Jacobi bracket $\{\,,\,\}_L$, which
is given by (\ref{19}) or (\ref{20}).
\end{theorem}

By applying Theorem \ref{th-dir-jacquot} to the case of the
generalized Lie bialgebroid defined by a Jacobi structure
$(\Lambda,E)$ on $M$ (Theorem \ref{Jac-bialg}) we deduce :

\begin{corol}\label{corol-Jac}
Let $(M,\Lambda,E)$ be a Jacobi manifold, $\big((TM\times
\R,[\,,\,],\pi, (0,1)),(T^{\ast}M\times \R, [\,,\,]_{(\Lambda,E)},
\pi\circ (\Lambda,E)^{\#}, (-E,0))\big)$ the associated
generalized Lie bialgebroid and $L$ a reducible Dirac structure
for the generalized Courant algebroid $\big((TM\times \R)\oplus
(T^{\ast}M\times \R), \lcf\,,\,\rcf, (\,,\,)_+, \pi + \pi \circ
(\Lambda,E)^{\#}, (0,1)+(-E,0)\big)$. We suppose that $1$ is an
$L$-admissible function. Then $L$ induces a Jacobi structure on
$M/\mathcal{F}$, where $\mathcal{F}$ is the foliation of $M$
defined by the distribution $\pi(D)$, $D=L\cap (TM\times \R)$,
which is exactly the Jacobi structure defined by $\{\,,\,\}_L$.
\end{corol}

\begin{remark}
{\rm In the context of Corollary \ref{corol-Jac}, the condition}
"$1$ is an $L$-admissible function" {\rm is equivalent to the one}
"$D$ has only sections of type $(X,0)$ with $X\in \Gamma(TM)$".
{\rm In fact, according to Lemma \ref{lemma-1-admis}, $1\in
C_L^{\infty}(M,\R)$ if and only if, for any $(X,f)\in \Gamma(D)$,
$\langle (0,1), (X,f)\rangle \stackrel{(\ref{eq1})}{=} 0
\Leftrightarrow f = 0$.}
\end{remark}

Taking into account Corollary \ref{corol-bialg-Jac}, Definition
\ref{Def-admis} and (\ref{20}), we can easily establish :

\begin{proposition}\label{prop-Jac-admis}
Under the assumptions of Corollary \ref{corol-Jac},
\begin{enumerate}
\item if $L = graph(\Lambda',E')^{\#}$ is the graph of a
$(TM\times \R)$-bivector field $(\Lambda',E')$ on $M$, then
$C_L^{\infty}(M,\R)=C^{\infty}(M,\R)$ and, for all $f,g \in
C_L^{\infty}(M,\R)$,
\begin{equation}\label{e1}
\{f,g\}_L = \{f,g\}_{(\Lambda',E')} + \{f,g\}_{(\Lambda,E)}\, ;
\end{equation}
\item if $L = D\oplus D^{\bot}$ is a null Dirac structure, then
$C_L^{\infty}(M,\R) = \{f \in C^{\infty}(M,\R) \,/\, (\delta f,f)
\in \Gamma(D^{\bot})\}$ and, for all $f,g \in C_L^{\infty}(M,\R)$,
\begin{equation}\label{e2}
\{f,g\}_L = \{f,g\}_{(\Lambda,E)}\, ;
\end{equation}
\item if $L = D\oplus graph(\Lambda',E')^{\#}\vert_{D^{\bot}}$ is
defined by a characteristic pair $(D, (\Lambda',E'))$, then
$C_L^{\infty}(M,\R) = \{f \in C^{\infty}(M,\R) \,/\, (\delta f,f)
\in \Gamma(D^{\bot})\}$ and, for all $f,g \in C_L^{\infty}(M,\R)$,
\begin{equation}\label{e3}
\{f,g\}_L = \{f,g\}_{(\Lambda',E')} + \{f,g\}_{(\Lambda,E)}.
\end{equation}
\end{enumerate}
\end{proposition}

In what follows, we will prove that in the context of {\it
"generalized Lie bialgebroids - Jacobi structures"}, as in the
context of {\it "Lie bialgebroids - Poisson structures"}
(\cite{lwx2}), the converse result of Corollary \ref{corol-Jac}
also holds.

\begin{theorem}\label{th-jacquot-Dir}
Let $(M,\Lambda,E)$ be a Jacobi manifold, $\mathcal{F}$ a simple
foliation of $M$ defined by a Lie subalgebroid $D\subset TM\times
\R$ that has only sections of type $(X,0)$ and
$(\Lambda_{M/\mathcal{F}},E_{M/\mathcal{F}})$ a Jacobi structure
on the quotient manifold $M/\mathcal{F}$. Then
$(M/\mathcal{F},\Lambda_{M/\mathcal{F}},E_{M/\mathcal{F}})$
defines a reducible Dirac structure $L$ in $(TM\times \R)\oplus
(T^{\ast}M\times \R)$ such that $L\cap (TM\times \R) =D$, $1\in
C_L^{\infty}(M,\R)$ and the Jacobi structure induced by $L$ on
$M/\mathcal{F}$, in the sense of Corollary \ref{corol-Jac}, is the
initially given $(\Lambda_{M/\mathcal{F}},E_{M/\mathcal{F}})$.
\end{theorem}
{\bf Proof :} We make the proof in several steps.

\vspace{2mm} \noindent {\it First step :} Let $D\subset TM\times
\R$ be a Lie subalgebroid of $(TM\times \R,[\,,\,],\pi)$, which
has only sections of type $(X,0)$, such that $\pi(D)$ defines a
simple foliation $\mathcal{F}$ of $M$ and let $D^{\bot}$ be its
conormal bundle :
\begin{eqnarray}\label{eq3}
D^{\bot} & = & \{(\alpha,g)\in T^{\ast}M\times
\R\,/\,\langle(\alpha,g),(X,0)\rangle = \langle \alpha,X\rangle =
0,\; \forall (X,0)\in D\}  \nonumber \\
& = & \pi(D)^{\bot}\times \R.
\end{eqnarray}
We suppose that the quotient manifold $M/\mathcal{F}$ is endowed
with a Jacobi structure
$(\Lambda_{M/\mathcal{F}},E_{M/\mathcal{F}})$ and we denote by $p
: M\to M/\mathcal{F}$ the canonical projection.

\vspace{2mm} \noindent {\it Second step :} We keep under control
the fact that $p : M \to M/\mathcal{F}$ is not a Jacobi map by
defining a "difference" bracket $\{\,,\,\}_1 :
C^{\infty}(M/\mathcal{F},\R)\times C^{\infty}(M/\mathcal{F},\R)
\to C^{\infty}(M,\R)$ as follows :
\begin{equation}\label{31}
\{f,g\}_1 = p^{\ast}\{f,g\}_{M/\mathcal{F}}-\{
p^{\ast}f,p^{\ast}g\}_{(\Lambda,E)}, \hspace{4mm}\forall \, f,g
\in C^{\infty}(M/\mathcal{F},\R).
\end{equation}
Obviously, $\{\,,\,\}_1$ is a bilinear, skew-symmetric, first
order differential operator on each of its arguments. Thus,
$\{\,,\,\}_1$ induces a skew-symmetric bilinear form
$(\Lambda_1,E_1)$ on $T^{\ast}(M/\mathcal{F})\times \R$ so that,
for all $f,g \in C^{\infty}(M/\mathcal{F},\R)$,
$$
\{f,g\}_1 = \Lambda_1(\delta f,\delta g) + \langle f\delta g -
g\delta f, E_1\rangle.
$$
In turn, $(\Lambda_1,E_1)$ induces a vector bundle map
$(\Lambda_1,E_1)^{\#} : T^{\ast}(M/\mathcal{F})\times \R \to
T(M/\mathcal{F})\times \R$. But, $T^{\ast}(M/\mathcal{F})\times \R
\cong \pi(D)^{\bot} \times \R \stackrel{(\ref{eq3})}{=}D^{\bot}$
and $T(M/\mathcal{F})\times \R \cong (TM/\pi(D))\times \R \cong
(TM\times \R)/D$. Consequently, we can consider that
$(\Lambda_1,E_1)^{\#} : D^{\bot} \to (TM\times \R)/D$.

\vspace{2mm} \noindent {\it Third step :} We denote by $pr :
TM\times \R \to (TM\times \R)/D$ the natural projection and we
define a subbundle $L\subset (TM\times \R)\oplus (T^{\ast}M\times
\R)$ by
\begin{equation}\label{32}
L = \{(X,f)+(\alpha,g)\in (TM\times \R)\oplus D^{\bot}\,/\,
pr(X,f) = (\Lambda_1,E_1)^{\#}(\alpha,g)\}.
\end{equation}
By construction, $L$ is maximally isotropic,
$C^{\infty}_L(M,\R)\cong C^{\infty}(M/\mathcal{F},\R)$ and $1\in
C^{\infty}_L(M,\R)$. Effectively, by a straightforward calculation
we show that, for any
$e_1=(X_1,f_1)+(\alpha_1,g_1),e_2=(X_2,f_2)+(\alpha_2,g_2)\in L$,
$(e_1,e_2)_+ = 0$ and $f\in C^{\infty}_L(M,\R)\Leftrightarrow
d^{(0,1)}f=(\delta f,f)\in \Gamma(D^{\bot})\cong
\Gamma(T^{\ast}(M/\mathcal{F})\times \R)\Leftrightarrow f\in
C^{\infty}(M/\mathcal{F},\R)$. Obviously, $1\in
C^{\infty}_L(M,\R)$ since $(0,1)\in \Gamma(D^{\bot})\cong
\Gamma(T^{\ast}(M/\mathcal{F})\times \R)$. Also, by Definition
\ref{Def-admis}, $f\in C^{\infty}_L(M,\R)$ if and only if there
exists $(Y_{f},\varphi_{f})\in \Gamma(TM\times \R)$ such that
$e_{f}=(Y_{f},\varphi_{f}) + (\delta f,f)\in \Gamma(L)$. Hence, we
have that $\Gamma(L)$ is spanned by all the sections of the type
$he_{f}$, where $h\in C^{\infty}(M,\R)$ and $f\in
C^{\infty}_L(M,\R)$. To verify the integrability of $L$, it
suffices to verify the closeness of the bracket $\lcf \,,\,\rcf$
for the sections of $L$ of the form $e_{f}=(Y_{f},\varphi_{f}) +
(\delta f,f)$ with $f\in C^{\infty}_L(M,\R)$, since, according to
(\ref{3def}) and because $L$ is isotropic,
\begin{eqnarray*}
\lcf e_{f},he_{g}\rcf & = & h\lcf e_{f},e_{g}\rcf +
(\rho(e_{f})h)e_{g} - (e_{f},e_{g})_+\mathcal{D}h  \\
 & = & h\lcf e_{f},e_{g}\rcf +
(\rho(e_{f})h)e_{g},
\end{eqnarray*}
for all $e_{f},e_{g}\in \Gamma(L)$, with $f,g\in
C_L^{\infty}(M,\R)$, and $h\in C^{\infty}(M,\R)$.

Let $f,g\in C_L^{\infty}(M,\R)$ be two $L$-admissible functions.
Since $C_L^{\infty}(M,\R)\cong C^{\infty}(M/\mathcal{F},\R)$,
$\{f,g\}_{M/\mathcal{F}}\in C_L^{\infty}(M,\R)$, i.e there is
$(Y_{\{f,g\}_{M/\mathcal{F}}},\varphi_{\{f,g\}_{M/\mathcal{F}}})\in
\Gamma(TM\times \R)$ such that $e_{\{f,g\}_{M/\mathcal{F}}}=
(Y_{\{f,g\}_{M/\mathcal{F}}},\varphi_{\{f,g\}_{M/\mathcal{F}}})$
$+ (\delta \{f,g\}_{M/\mathcal{F}},\{f,g\}_{M/\mathcal{F}})\in
\Gamma(L)$. We show that
\begin{equation}\label{35}
\{f,g\}_{M/\mathcal{F}} = \rho^{\theta}(e_{f})g
\stackrel{(\ref{19})}{=} : \{f,g\}_L.
\end{equation}
Effectively,
\begin{eqnarray*}
\{f,g\}_L = \rho^{\theta}(e_{f})g & = & \big[\big(\pi^{(0,1)} +
(\pi \circ
(\Lambda,E)^{\#})^{(-E,0)}\big)\big((Y_{f},\varphi_{f})+(\delta
f,f)\big)\big]g  \\
& = & \big(Y_{f}+\varphi_{f}+\Lambda^{\#}(\delta f)+ f E-\langle
\delta f, E\rangle \big)g  \\
& = & (pr(Y_{f},\varphi_{f}) + the\; component\; of\;
(Y_{f},\varphi_{f})\; on \; D)g + \{f,g\}_{(\Lambda,E)}  \\
& = & \langle (\delta g,g), (\Lambda_1,E_1)^{\#}(\delta
f,f)\rangle +
\{f,g\}_{(\Lambda,E)}  \\
& = & \{f,g\}_1 + \{f,g\}_{(\Lambda,E)}  \\
& \stackrel{(\ref{31})}{=}& \{f,g\}_{M/\mathcal{F}}.
\end{eqnarray*}
On the other hand, since $\{\,,\,\}_{M/\mathcal{F}}$ is a Jacobi
bracket, thus it verifies the Jacobi identity, for any $f,g,h \in
C_L^{\infty}(M,\R)\cong C^{\infty}(M/\mathcal{F},\R)$,
\begin{eqnarray}\label{36}
\rho^{\theta}(\lcf e_{f},e_{g}\rcf - e_{\{f,g\}_{M/\mathcal{F}}})h
& = & \rho^{\theta}(\lcf e_{f},e_{g}\rcf)h -
\rho^{\theta}(e_{\{f,g\}_{M/\mathcal{F}}})h  \nonumber \\
& \stackrel{(\ref{33})}{=}&
[\rho^{\theta}(e_{f}),\rho^{\theta}(e_{g})]h-
\rho^{\theta}(e_{\{f,g\}_{M/\mathcal{F}}})h  \nonumber \\
& = & \rho^{\theta}(e_{f})(\rho^{\theta}(e_{g})h)-
\rho^{\theta}(e_{g})(\rho^{\theta}(e_{f})h)-
\rho^{\theta}(e_{\{f,g\}_{M/\mathcal{F}}})h  \nonumber \\
& \stackrel{(\ref{35})}{=}&
\{f,\{g,h\}_{M/\mathcal{F}}\}_{M/\mathcal{F}} -
\{g,\{f,h\}_{M/\mathcal{F}}\}_{M/\mathcal{F}} - \{\{f,g\}_{M/\mathcal{F}},h\}_{M/\mathcal{F}}  \nonumber \\
& = & 0.
\end{eqnarray}
From the proof of Proposition \ref{pr-admis}, we have that the
component of $\lcf e_{f},e_{g}\rcf$ in $\Gamma(T^{\ast}M \times
\R)$ is $d^{(0,1)}\{f,g\}_{M/\mathcal{F}}$, therefore $\lcf
e_{f},e_{g}\rcf - e_{\{f,g\}_{M/\mathcal{F}}} \in \Gamma(TM\times
\R)$. So, (\ref{36}) means that $\rho^{\theta}(\lcf
e_{f},e_{g}\rcf - e_{\{f,g\}_{M/\mathcal{F}}})\in \Gamma(D)$. But
$\rho^{\theta}(\lcf e_{f},e_{g}\rcf - e_{\{f,g\}_{M/\mathcal{F}}})
= \pi^{(0,1)}(\lcf e_{f},e_{g}\rcf - e_{\{f,g\}_{M/\mathcal{F}}})=
\lcf e_{f},e_{g}\rcf - e_{\{f,g\}_{M/\mathcal{F}}}$ and $\Gamma(D)
\subset \Gamma(L)$. Consequently, $\lcf e_{f},e_{g}\rcf -
e_{\{f,g\}_{M/\mathcal{F}}} \in \Gamma(L)$ which implies $\lcf
e_{f},e_{g}\rcf \in \Gamma(L)$, whence the integrability of $L$.

For the constructed $L$, $L\cap (TM\times \R) = \{(X,f)+(0,0) \in
(TM\times \R)\oplus \{(0,0)\} \, /\,
pr(X,f)=(\Lambda_1,E_1)^{\#}(0,0)\}=\{(X,f) \in TM\times \R \, /\,
pr(X,f)=(0,0)\}=D$ and the induced Jacobi structure on
$M/\mathcal{F}$, in the sense of Corollary \ref{corol-Jac}, is the
initially given $(\Lambda_{M/\mathcal{F}},E_{M/\mathcal{F}})$
(see, (\ref{35})). $\bullet$

\begin{remark}\label{R2}
{\rm The condition "$D$ {\it has only sections of type $(X,0)$}"
is indispensable in order to ensure that the constant function $1$
is an $L$-admissible function for the constructed $L$. In the
opposite case, i.e if $D$ has at least one section of type $(X,f)$
with $f\neq 0$, we will have that there exists at least one
section of $D$, the section $(X,f)$, such that $\langle
(0,1),(X,f)\rangle = f \neq 0$ and, according to Lemma
\ref{lemma-1-admis}, this implies that $1$ is not an
$L$-admissible function. Hence, $(C_L^{\infty}(M,\R),
\{\,,\,\}_L)$ can not be a Jacobi algebra and $C_L^{\infty}(M,\R)$
does not coincide with $C^{\infty}(M/\mathcal{F},\R)$. Thus, for a
Lie subalgebroid $D$ of $TM\times \R$ that has at least one
section $(X,f)$ with $f\neq 0$ we can not construct a reducible
Dirac subbundle $L\subset (TM\times \R)\oplus (T^{\ast}M\times
\R)$ which induces, in the sense of Corollary \ref{corol-Jac}, a
Jacobi structure on $M/\mathcal{F}$.}
\end{remark}

In conclusion, we have proved :

\begin{theorem}\label{basic-th}
Let $(M,\Lambda,E)$ be a Jacobi manifold. There is a one-one
correspondence between reducible Dirac subbundles of the
generalized Courant algebroid $\big((TM\times \R)\oplus
(T^{\ast}M\times \R), \lcf\,,\,\rcf, (\,,\,)_+, \pi + \pi \circ
(\Lambda,E)^{\#}, (0,1)+(-E,0)\big)$ for which $1$ is an
admissible function and quotient Jacobi manifolds $M/\mathcal{F}$
of $M$, where $\mathcal{F}$ is a simple foliation of $M$ defined
by a Lie subalgebroid $D\subset TM\times \R$ that has sections
only of type $(X,0)$.
\end{theorem}

\begin{remark}\label{R3}
{\rm If, in the proof of Theorem \ref{th-jacquot-Dir}, $p :
(M,\Lambda,E) \to
(M/\mathcal{F},\Lambda_{M/\mathcal{F}},E_{M/\mathcal{F}})$ is a
Jacobi map, then $(\Lambda_1,E_1)=(0,0)$. Hence, in this case, $L
= D\oplus D^{\bot}$ is a null Dirac structure. Thus :}
\end{remark}

\begin{corol}\label{C1}
A Lie subalgebroid $D\subset TM\times \R$ which has only sections
of type $(X,0)$ defines a simple foliation $\mathcal{F}$ of
$(M,\Lambda,E)$ such that $p : (M,\Lambda,E) \to
(M/\mathcal{F},\Lambda_{M/\mathcal{F}},E_{M/\mathcal{F}})$ is a
Jacobi map if and only if $L = D\oplus D^{\bot}$.
\end{corol}

\begin{remark}\label{R4}
{\rm In the case where $D=\{(0,0)\}$, a Jacobi structure on
$M/\mathcal{F} \cong M$ is a new Jacobi structure $(\Lambda',E')$
on $M$ and the constructed $L$ is the graph of
$(\Lambda'-\Lambda,E'-E)$. Since, by construction, $L$ is a Dirac
subbundle of $(TM\times \R)\oplus (T^{\ast}M\times \R)$,
$(\Lambda'-\Lambda,E'-E)$ is a Jacobi structure on $M$
(\cite{jj}), fact which implies that $(\Lambda,E)$ and
$(\Lambda',E')$ are compatible Jacobi structures in the sense of
\cite{j1}}.
\end{remark}

\vspace{3mm} \noindent {\bf A geometric interpretation of
Corollary \ref{corol-charpair} :} In the context of this
paragraph, Corollary \ref{corol-charpair} can be formulated as :
{\it Let $(M,\Lambda,E)$ be a Jacobi manifold, $\big ((TM\times
\R, (0,1)),(T^{\ast}M\times \R,(-E,0)),(\Lambda,E)\big )$ the
associated triangular generalized Lie bialgebroid over $M$ and
$(\Lambda',E')$ a $(TM\times \R)$-bivector field such that
$L=D\oplus graph((\Lambda',E')^{\#}\vert_{D^{\bot}})$ is a maximal
isotropic subbundle of $(TM\times \R)\oplus (T^{\ast}M\times \R)$
with fixed characteristic pair $(D,(\Lambda',E'))$. Then $L$ is a
Dirac structure for $((TM\times \R)\oplus (T^{\ast}M\times
\R),(0,1)+(-E,0))$ if and only if
\begin{enumerate}
\item[(i)] $D$ is a Lie subalgebroid of $TM\times \R$ ;
\item[(ii)] $[(\Lambda +\Lambda',E+E'),(\Lambda
+\Lambda',E+E')]^{(0,1)}\equiv 0(mod D)$ ; \item[(iii)] for any
$(X,f)\in \Gamma(D)$,
$\mathcal{L}_{(X,f)}^{(0,1)}(\Lambda+\Lambda',E+E')\equiv 0(mod
D)$.
\end{enumerate}}
\noindent If $L=D\oplus graph((\Lambda',E')^{\#}\vert_{D^{\bot}})$
is a reducible Dirac structure and $1$ is an $L$-admissible
function, after the proofs of Theorems \ref{th-dir-jacquot} and
\ref{th-jacquot-Dir}, we get that condition (iii) is equivalent to
that $(\Lambda+\Lambda',E+E')$ can be reduced to a $(TM\times
\R)/D \cong (T(M/\mathcal{F})\times \R)$-bivector field on
$M/\mathcal{F}$ and the condition (ii) is equivalent to the fact
that the reduced bivector field is a Jacobi structure on
$M/\mathcal{F}$. Furthermore, by Proposition \ref{prop-Jac-admis}
(case 3) we get that the induced Jacobi structure on
$M/\mathcal{F}$ is exactly the one defined by the bracket of
$L$-admissible functions. Consequently, it is the Jacobi structure
induced on $M/\mathcal{F}$ by $L$ in the sense of Corollary
\ref{corol-Jac}.

\section{Dirac structures and Jacobi reduction}
In this paragraph, we will establish a Jacobi reduction theorem in
terms of Dirac structures. For its proof, we need to adapt the
results concerning the pull-back Dirac structures of a Lie
bialgebroid (\cite{liu}) to the pull-back Dirac structures for a
generalized Lie bialgebroid.

\begin{proposition}\label{pull-back}
Let $(A_1,\phi_1)$ be a Lie algebroid over a differentiable
manifold $M_1$ with an 1-cocycle,
$((A_2,\phi_2),(A_2^{\ast},W_2),P_2)$ a triangular generalized Lie
bialgebroid over a differentiable manifold $M_2$ and $\Phi : A_1
\to A_2$ a Lie algebroid morphism of constant rank, which covers a
surjective map between the bases, such that
$\Phi^{\ast}(\phi_2)=\phi_1$. Then the following two statements
are equivalent.
\begin{enumerate}
\item There exists a Dirac structure for the triangular
generalized Lie bialgebroid $((A_1,\phi_1),(A_1^{\ast},0),0)$
whose characteristic pair is $(\ker \Phi, P_1)$ and $\Phi(P_1) =
P_2$. \item $\mathrm{Im}P_2^{\#} \subseteq \mathrm{Im}\Phi$.
\end{enumerate}
We note that, since $\Phi : A_1 \to A_2$ is a Lie algebroid
morphism such that $\Phi^{\ast}(\phi_2)=\phi_1$ then, for any
$P\in \Gamma(\bigwedge^pA_1)$ and $Q\in \Gamma(\bigwedge^qA_1)$,
$\Phi([P,Q]^{\phi_1}) = [\Phi(P),\Phi(Q)]^{\phi_2}$.
\end{proposition}
{\bf Proof :} According to Corollary \ref{corol-charpair}, it
suffices to show that the following two statements are equivalent.
\begin{enumerate}
\item There exists $P_1\in \Gamma(\bigwedge^2A_1)$ such that
$\Phi(P_1) = P_2$ and
\begin{enumerate}
\item $\ker \Phi$ is a Lie subalgebroid of $A_1$ ; \item
$[0+P_1,0+P_1]^{\phi_1}\equiv 0 (mod \ker\Phi)\Leftrightarrow
[P_1,P_1]^{\phi_1}\equiv 0 (mod \ker\Phi)$ ; \item for any $X\in
\Gamma(\ker \Phi)$, $\mathcal{L}_X^{\phi_1}(0+P_1)\equiv 0 (mod
\ker \Phi)\Leftrightarrow \mathcal{L}_X^{\phi_1}(P_1)\equiv 0 (mod
\ker \Phi)$.
\end{enumerate}
\item $\mathrm{Im}P_2^{\#} \subseteq \mathrm{Im}\Phi$.
\end{enumerate}
Obviously, $\ker \Phi$ is a Lie subalgebroid of $A_1$ since, for
all $X,Y\in \Gamma(\ker \Phi)$, $\Phi ([X,Y])
=[\Phi(X),\Phi(Y)]=[0,0]=0$, which means that $[X,Y]\in
\Gamma(\ker \Phi)$. On the other hand, the subbundle $\ker
\Phi^{\bot} = \{\alpha \in A_1^{\ast}\,/\, \langle \alpha,X\rangle
= 0,\; \forall \, X\in \ker \Phi\}$ of $A_1^{\ast}$ can be
identified with the dual bundle $(A_1/\ker \Phi)^{\ast}$ of
$A_1/\ker \Phi$. Also, $\ker \Phi^{\bot} =
\mathrm{Im}\Phi^{\ast}$, where $\Phi^{\ast} : A_2^{\ast}\to
A_1^{\ast}$ is the dual map of $\Phi$. Effectively, it is clear
that, $\mathrm{Im}\Phi^{\ast}\subseteq \ker \Phi^{\bot}$ and,
since $\Phi$ is of constant rank, $\dim \mathrm{Im}\Phi^{\ast}
=\dim \ker \Phi^{\bot}$, thus $\mathrm{Im}\Phi^{\ast} =\ker
\Phi^{\bot}\cong (A_1/\ker \Phi)^{\ast}$. Hence, $\Phi^{\ast} :
A_2^{\ast} \to (A_1/\ker \Phi)^{\ast}$ is a surjective map, i.e.,
for any $\bar{\alpha}_1,\bar{\beta}_1 \in \Gamma((A_1/\ker
\Phi)^{\ast})$, there exist $\alpha_2,\beta_2 \in
\Gamma(A_2^{\ast})$ such that
$\bar{\alpha}_1=\Phi^{\ast}(\alpha_2)$ and
$\bar{\beta}_1=\Phi^{\ast}(\beta_2)$. If there is some
$\bar{P}_1\in \Gamma(\bigwedge^2(A_1/\ker \Phi))$ which is
$\Phi$-related to $P_2$, i.e. $\Phi(\bar{P}_1)=P_2$, then it
should be defined, for all $\bar{\alpha}_1,\bar{\beta}_1 \in
\Gamma((A_1/\ker \Phi)^{\ast})$, by
$$
\bar{P}_1(\bar{\alpha}_1,\bar{\beta}_1) = P_2(\alpha_2,\beta_2).
$$
It is clear that $\bar{P}_1$ is well-defined if and only if $\ker
\Phi^{\ast} \subseteq \ker P_2^{\#}$, or equivalently, if and only
if $\mathrm{Im}P_2^{\#} \subseteq \mathrm{Im}\Phi$. Let $P_1$ be
an arbitrary representative of $\bar{P}_1$ in
$\Gamma(\bigwedge^2A_1)$. Since $\Phi : A_1\to A_2$ is a Lie
algebroid morphism such that $\Phi^{\ast}(\phi_2)=\phi_1$ and
$((A_2,\phi_2),(A_2^{\ast},W_2),P_2)$ is a triangular generalized
Lie bialgebroid, we have that
$$
\Phi([P_1,P_1]^{\phi_1})=[\Phi(P_1),\Phi(P_1)]^{\phi_2}=[P_2,P_2]^{\phi_2}
= 0 \;\Leftrightarrow \; [P_1,P_1]^{\phi_1}\equiv 0(mod \ker
\Phi).
$$
Moreover, for any $X\in \Gamma(\ker\Phi)$,
$$
\Phi(\mathcal{L}_X^{\phi_1}P_1)=\Phi([X,P_1]^{\phi_1})=
[\Phi(X),\Phi(P_1)]^{\phi_2}=[0,\Phi(P_1)]^{\phi_2}=0 \;
\Leftrightarrow \; \mathcal{L}_X^{\phi_1}P_1\equiv 0 (mod
\ker\Phi).
$$
Consequently, there exists $P_1\in \Gamma(\bigwedge^2A_1)$ such
that $\Phi(P_1) = P_2$ and $(\ker\Phi,P_1)$ defines a Dirac
structure for the triangular generalized Lie bialgebroid
$((A_1,\phi_1),(A_1^{\ast},0),0)$ if and only if
$\mathrm{Im}P_2^{\#} \subseteq \mathrm{Im}\Phi$. $\bullet$

\vspace{3mm} \noindent {\bf Reduction of Jacobi manifolds :} Let
$(M,\Lambda,E)$ be a Jacobi manifold, $N\subseteq M$ a submanifold
of $M$ and $i : N\hookrightarrow M$ the canonical inclusion,
$D\subset TM\times \R$ a Lie subalgebroid of $(TM\times \R,
[\,,\,],\pi)$ that has only sections of type $(X,0)$ and $D_0 =
D\cap(TN\times \R)$. We suppose that $D$ and $D_0$ define,
respectively, a simple foliation $\mathcal{F}$ of $M$ and a simple
foliation $\mathcal{F}_0$ of $N$ and we denote by $p : M\to
M/\mathcal{F}$ and $p_0 : N \to N/\mathcal{F}_0$ the canonical
projections. Thus, we have the following commutative diagram :
\begin{equation}\label{diagram1}
\begin{array}{rcl}
 N & \stackrel{i}{\hookrightarrow} & M  \\
 & & \\
p_0  \downarrow &  & \downarrow  p \\
&  & \\
N/\mathcal{F}_0 & \stackrel{\varphi}{\rightarrow} & M/\mathcal{F}
\end{array}
\end{equation}
Since any leaf of $\mathcal{F}_0$ is a connected component of the
intersection between $N$ and some leaf of $\mathcal{F}$, we can
always suppose, under some clean intersection condition, that
$\varphi : N/\mathcal{F}_0 \to M/\mathcal{F}$ is an immersion,
locally injective.

We consider $L=D\oplus D^{\bot}$ and we suppose that $L$ is a null
Dirac structure for the triangular generalized Lie bialgebroid
$\big ((TM\times \R, (0,1)),(T^{\ast}M\times
\R,(-E,0)),(\Lambda,E)\big )$. By the hypothesis on $D$, we have
that $L$ is also reducible and that $1$ is an $L$-admissible
function. Then, by Corollary \ref{corol-Jac}, we get that $L$
induces a Jacobi structure
$(\Lambda_{M/\mathcal{F}},E_{M/\mathcal{F}})$ on $M/\mathcal{F}$
and by Corollary \ref{C1} and Remark \ref{R3}, we obtain that $p :
(M,\Lambda,E)\to
(M/\mathcal{F},\Lambda_{M/\mathcal{F}},E_{M/\mathcal{F}})$ is a
Jacobi map. We consider the triangular generalized Lie
bialgebroids $\big((T(M/\mathcal{F})\times \R, (0,1)),
(T^{\ast}(M/\mathcal{F})\times \R,(-E_{M/\mathcal{F}},0)),
(\Lambda_{M/\mathcal{F}},E_{M/\mathcal{F}})\big)$ over
$M/\mathcal{F}$ and $\big((TN\times \R, (0,1)), (T^{\ast}N\times
\R,(0,0)), (0,0)\big)$ over $N$. We note that any function
$\bar{f}\in C^{\infty}(N,\R)$ can be seen as the image by $(p
\circ i)^{\ast}$ of a function $f\in
C^{\infty}(M/\mathcal{F},\R)$, i.e $\bar{f} =(p \circ i)^{\ast}f$.
Since $\mathcal{F}$ is a regular foliation, $p$ has constant rank,
thus the map $p \circ i : N \to M/\mathcal{F}$ has also constant
rank. Hence, the application $\Phi : TN\times \R \to
T(M/\mathcal{F})\times \R \cong (TM\times \R)/D$ defined, for any
$(X,\bar{f})\in \Gamma(TN\times \R)$, $\bar{f} = (p\circ
i)^{\ast}f$ with $f\in C^{\infty}(M/\mathcal{F},\R)$, by
\begin{equation}\label{37}
\Phi(X,\bar{f})=((p \circ i)_{\ast}X,f),
\end{equation}
can be considered as a Lie algebroid morphism of constant rank
such that $\Phi^{\ast}(0,1)=(0,1)$ and $\ker \Phi = D\cap
(TN\times \R)=D_0$. Therefore, by Proposition \ref{pull-back},
there exists a pull-back Dirac structure $L_0$ for the triangular
generalized Lie bialgebroid $\big((TN\times \R, (0,1)),
(T^{\ast}N\times \R,(0,0)), (0,0)\big)$ with characteristic pair
$(D_0,(\Lambda_N,E_N))$ satisfying $\Phi(\Lambda_N,E_N) =
(\Lambda_{M/\mathcal{F}},E_{M/\mathcal{F}})$ if and only if
$\mathrm{Im}(\Lambda_{M/\mathcal{F}},E_{M/\mathcal{F}})^{\#}\subseteq
\mathrm{Im}\Phi$ holds on $T(M/\mathcal{F})\times \R$, i.e.
\begin{equation}\label{38}
\Gamma((\Lambda_{M/\mathcal{F}},E_{M/\mathcal{F}})^{\#}(D^{\bot}))
\subseteq \{((p\circ i)_{\ast}X,f) \,/\, X \in \Gamma(TN)\;
\mathrm{and}\; f\in C^{\infty}(M/\mathcal{F},\R)\}.
\end{equation}
But, $p : (M,\Lambda,E) \to
(M/\mathcal{F},\Lambda_{M/\mathcal{F}},E_{M/\mathcal{F}})$ being a
Jacobi map, $(\Lambda_{M/\mathcal{F}},E_{M/\mathcal{F}}) =
p_{\ast}(\Lambda,E)$. Thus, on the submanifold $N\subseteq M$, by
identifying $i_{\ast}(TN)$ with $TN$, condition (\ref{38}) is
equivalent to
\begin{equation}\label{43}
(\Lambda,E)^{\#}(D^{\bot})\subseteq TN\times \R + D.
\end{equation}
On the other hand, since $D=\pi(D)\times \{0\}$, $D^{\bot} =
\pi(D)^{\bot}\times \R$, consequently, (\ref{43}) is equivalent to
\begin{equation}\label{eq4}
\Lambda^{\#}(\pi(D)^{\bot})\subseteq TN + \pi(D) \hspace{3mm}
\mathrm{and} \hspace{3mm} E\vert_N \in \Gamma(TN + \pi(D)).
\end{equation}
Also, since $L_0 = D_0 \oplus graph
(\Lambda_N,E_N)^{\#}\vert_{D_0^{\bot}}$ is a reducible Dirac
structure of $((TN\times \R)\oplus (T^{\ast}N\times \R),
(0,1)+(0,0))$ and $1\in C_{L_0}^{\infty}(N,\R)$, it induces a
Jacobi structure $(\Lambda_{N/\mathcal{F}_0},E_{N/\mathcal{F}_0})$
on $N/\mathcal{F}_0$ (see, Corollary \ref{corol-Jac}) such that
$(\Lambda_{N/\mathcal{F}_0},E_{N/\mathcal{F}_0}) =
p_{0\ast}(\Lambda_N,E_N)$ (see, Corollary \ref{corol-charpair} and
its geometric interpretation). By the above results and by the
commutativity of the diagram (\ref{diagram1}), we obtain :
\begin{eqnarray}\label{41}
(\Lambda_{M/\mathcal{F}},E_{M/\mathcal{F}}) & = & \Phi(\Lambda_N,E_N)  \nonumber \\
                & = & \big((p \circ i)_{\ast}\Lambda_N,(p \circ i)_{\ast}E_N\big) \nonumber \\
                & = & \big((\varphi \circ p_{0})_{\ast}\Lambda_N,
                      (\varphi \circ p_{0})_{\ast}E_N\big) \nonumber \\
                & = &
                \big(\varphi_{\ast}(p_{0\ast}\Lambda_N),\varphi_{\ast}(p_{0\ast}E_N)\big) \nonumber
                \\
                & = & \big(\varphi_{\ast}\Lambda_{N/\mathcal{F}_0},
                \varphi_{\ast}E_{N/\mathcal{F}_0}\big) \nonumber
                \\
                & = &
                \varphi_{\ast}(\Lambda_{N/\mathcal{F}_0},E_{N/\mathcal{F}_0}),
\end{eqnarray}
which means that $\varphi :
(N/\mathcal{F}_0,\Lambda_{N/\mathcal{F}_0},E_{N/\mathcal{F}_0})
\to (M/\mathcal{F}, \Lambda_{M/\mathcal{F}}, E_{M/\mathcal{F}})$
is a Jacobi map.

\vspace{1mm} The above study led us to the following theorem :

\begin{theorem}[{\bf Reduction Theorem of Jacobi manifolds}]\label{th-reduction}
Let $(M,\Lambda,E)$ be a Jacobi manifold, $N\subseteq M$ a
submanifold of $M$, $D\subset TM\times \R$ a Lie subalgebroid of
$(TM\times \R, [\,,\,],\pi)$ that has only sections of type
$(X,0)$ and $D_0 = D\cap(TN\times \R)$. We suppose that $D$ and
$D_0$ define, respectively, a simple foliation $\mathcal{F}$ of
$M$ and a simple foliation $\mathcal{F}_0$ of $N$ and that $L =
D\oplus D^{\bot}$ is a reducible Dirac structure for the
triangular generalized Lie bialgebroid $\big ((TM\times \R,
(0,1)),(T^{\ast}M\times \R,(-E,0)),(\Lambda,E)\big )$. Then, the
following two statements are equivalent.
\begin{enumerate}
\item There exists a Jacobi structure
$(\Lambda_{N/\mathcal{F}_0},E_{N/\mathcal{F}_0})$ on
$N/\mathcal{F}_0$ such that
$$
p_{\ast}(\Lambda,E) =
\varphi_{\ast}(\Lambda_{N/\mathcal{F}_0},E_{N/\mathcal{F}_0}).
$$
\item $\Lambda^{\#}(\pi(D)^{\bot})\subseteq TN + \pi(D)$ holds on
$N$ and $E\vert_N \in \Gamma(TN + \pi(D))$.
\end{enumerate}
\end{theorem}

\begin{remarks}\label{Rs1}
\end{remarks}
\vspace{-2mm}\noindent {\bf 1.} We remark that, in the context of
the Reduction Theorem \ref{th-reduction}, the initial Jacobi
manifold $(M,\Lambda,E)$ and the reduced Jacobi manifold
$(N/\mathcal{F}_0,\Lambda_{N/\mathcal{F}_0},E_{N/\mathcal{F}_0})$
are connected by means of the Jacobi manifold
$(M/\mathcal{F},\Lambda_{M/\mathcal{F}},E_{M/\mathcal{F}})$ with
two Jacobi maps.

\vspace{2mm} \noindent {\bf 2.} Reduction Theorem
\ref{th-reduction} holds for any reducible Dirac structure
$L\subset (TM\times \R)\oplus (T^{\ast}M\times \R)$ having a
characteristic pair $(D,(\Lambda',E'))$, i.e. $L = D \oplus graph
((\Lambda',E')^{\#}\vert_{D^{\bot}})$, such that $D$ has only
sections of type $(X,0)$, so $1\in C_L^{\infty}(M,\R)$.
Effectively, by Corollary \ref{corol-Jac} we get that $L$ induces
a Jacobi structure $(\Lambda_{M/\mathcal{F}},E_{M/\mathcal{F}})$
on $M/\mathcal{F}$ which is exactly the induced Jacobi structure
by $(\Lambda +\Lambda',E+E')$ (see, the geometric interpretation
of Corollary \ref{corol-charpair}). If
$(\Lambda_{M/\mathcal{F}},E_{M/\mathcal{F}})$ verifies (\ref{38})
or, equivalently, $(\Lambda +\Lambda', E+E')$ verifies
(\ref{eq4}), then, by Proposition \ref{pull-back}, there exists a
pull-back Dirac structure $L_0$ for $\big((TN\times \R,
(0,1)),(T^{\ast}N\times \R,(0,0)), (0,0)\big)$ with characteristic
pair $(D_0,(\Lambda_N,E_N))$ such that
$\Phi(\Lambda_N,E_N)=(\Lambda_{M/\mathcal{F}},E_{M/\mathcal{F}})$.
The reducible Dirac subbundle $L_0\subset (TN\times \R)\oplus
(T^{\ast}N\times \R)$ induces a Jacobi structure
$(\Lambda_{N/\mathcal{F}_0},E_{N/\mathcal{F}_0})$ on
$N/\mathcal{F}_0$ and
$$
p_{0\ast}(\Lambda_N,E_N) = (\Lambda_{N/\mathcal{F}_0},
E_{N/\mathcal{F}_0}).
$$
Applying the calculus of (\ref{41}) to the relation
$(\Lambda_{M/\mathcal{F}},E_{M/\mathcal{F}}) =
\Phi(\Lambda_N,E_N)$, we conclude that $\varphi : N/\mathcal{F}_0
\to M/\mathcal{F}$ is always a Jacobi map. But, the projection $p
: M \to M/\mathcal{F}$ is a Jacobi map if and only if $L$ is a
null Dirac structure, fact which is equivalent to
$(\Lambda',E')\equiv 0(mod D)$.

\vspace{2mm} \noindent {\bf 3.} As we have mentioned in
introduction, there are already several works treating the Jacobi
reduction problem. These results are, grosso-modo, equivalent to
the ones established by the second author in \cite{j2} and,
independently, by K. Mikami in \cite{mk}. They establish a
geometric Reduction Theorem for Jacobi manifolds by extending the
previous one proved by Marsden and Ratiu for Poisson manifolds
\cite{mr}, without mentioning Dirac structures. Precisely, they
prove :
\begin{theorem}\label{reduction-joana}
Let $(M,\Lambda,E)$ be a Jacobi manifold, $N$ a submanifold of $M$
and $\Delta$ a vector subbundle of $T_NM$ such that : (i) $\Delta
\cap TN$ defines a simple foliation $\mathcal{F}_0$ of $N$; (ii)
for any $f,g \in C^{\infty}(M,\R)$ with differentials $\delta f$
and $\delta g$, restricted to $N$, vanishing on $\Delta$, the
differential $\delta \{f,g\}_{(\Lambda,E)}$, restricted to $N$,
vanishes on $N$. Then, $(\Lambda,E)$ induces a unique Jacobi
structure $(\Lambda_{N/\mathcal{F}_0}, E_{N/\mathcal{F}_0})$ on
$N/\mathcal{F}_0$ if and only if
$\Lambda^{\#}(\Delta^{\bot})\subseteq TN + \Delta$ holds on $N$
and $E\vert_N \in \Gamma(TN + \Delta)$. The associated bracket of
$(\Lambda_{N/\mathcal{F}_0}, E_{N/\mathcal{F}_0})$ on
$C^{\infty}(N/\mathcal{F}_0,\R)$ is given, for any $f_0,g_0 \in
C^{\infty}(N/\mathcal{F}_0,\R)$ and any differentiable extensions
$f$ of $f_0 \circ p_0$ and $g$ of $g_0 \circ p_0$ with
differentials $\delta f$ and $\delta g$, restricted to $N$, vanish
on $\Delta$, by $\{f_0,g_0\}_{(\Lambda_{N/\mathcal{F}_0},
E_{N/\mathcal{F}_0})} \circ p_0 = \{f,g\}_{(\Lambda,E)}\circ i$,
where $p_0 : N \to N/\mathcal{F}_0$ is the canonical projection
and $i : N \to M$ is the canonical inclusion of $N$ into $M$.
\end{theorem}
We remark that the above Theorem is slightly different from
Theorem \ref{th-reduction}. In Theorem \ref{th-reduction} we
suppose that we have two simple foliations, a foliation
$\mathcal{F}$ of the initial phase space $M$ determined by
$\pi(D)$ and a foliation $\mathcal{F}_0$ of the considered
submanifold $N$ of $M$ determined by $\pi(D_0)=\pi(D)\cap TN$,
while in Theorem \ref{reduction-joana} we only suppose that we
have a subbundle $\Delta$ of $T_NM$ such that $\Delta \cap TN$
defines a simple foliation of $N$, also denoted by
$\mathcal{F}_0$. But, in both Theorems, the reducibility condition
$$
\Lambda^{\#}(\pi(D)^{\bot})\subseteq TN + \pi(D) \;\;\;
\mathrm{holds} \;\;  \mathrm{on}\;\; N \;\;\; \mathrm{and} \;\;\;
E\vert_N \in \Gamma(TN + \pi(D))
$$
is exactly the same. Thus, it is natural to ask : \emph{What is
the advantage of using Dirac structures in the study of Jacobi
reduction problem ?} The answer can be founded in Remarks 1 and 2
of this paragraph. By using reducible Dirac structures in this
study, we establish the existence, not only, of a reduced Jacobi
manifold $(N/\mathcal{F}_0,\Lambda_{N/\mathcal{F}_0},
E_{N/\mathcal{F}_0})$, but also of a quotient Jacobi manifold
$(M/\mathcal{F}_0,\Lambda_{M/\mathcal{F}_0}, E_{M/\mathcal{F}_0})$
which is always related with
$(N/\mathcal{F}_0,\Lambda_{N/\mathcal{F}_0}, E_{N/\mathcal{F}_0})$
by means of a Jacobi map ; very important fact when we treat
reduction problems. On the other hand, this study, in this
framework, allows us to investigate, in a future paper, the
\emph{Dirac reduction problem} and its relation with the one of
Jacobi, Poisson and symplectic structures.

\section{Applications and Examples}
{\bf 1. Jacobi submanifolds :} From Theorem \ref{th-reduction} we
obtain sufficient conditions under which a Jacobi structure
$(\Lambda,E)$ on a differentiable manifold $M$ induces a Jacobi
structure on a submanifold $N$ of $M$. Effectively, under the
assumptions of the above mentioned theorem, if $D_0 = D\cap (TN
\times \R) = \{(0,0)\}$ and $(\Lambda,E)^{\#}(D^{\bot})\subseteq
TN\times \R + D$ holds on $N$, then there exists a $(TN \times
\R)$-bivector field $(\Lambda_N,E_N)$ on $N$ such that $L_0 = D_0
\oplus graph(\Lambda_N,E_N)^{\#}\vert_{D_0^{\bot}} = graph
(\Lambda_N,E_N)^{\#}$ is a reducible Dirac structure for the
triangular generalized Lie bialgebroid $\big((TN\times
\R,(0,1)),(T^{\ast}N\times \R,(0,0)),0\big)$ and $\Phi
(\Lambda_N,E_N)=(\Lambda_{M/\mathcal{F}},E_{M/\mathcal{F}})$. But,
the fact {\it "$L_0 = graph(\Lambda_N,E_N)^{\#}$ is Dirac for
$\big((TN\times \R,(0,1)),(T^{\ast}N\times \R,(0,0)),0\big)$"} is
equivalent to the fact {\it "$(\Lambda_N,E_N)$ is a Jacobi
structure on $N$"} (see Proposition 5.2 in \cite{jj}) and
\begin{eqnarray*}
(\Lambda_{M/\mathcal{F}},E_{M/\mathcal{F}}) = \Phi (\Lambda_N,E_N)
& \Leftrightarrow & p_{\ast}(\Lambda,E) = (p \circ
i)_{\ast}(\Lambda_N,E_N) \\ & \Leftrightarrow &
p_{\ast}((\Lambda,E) - i_{\ast}(\Lambda_N,E_N)) = (0,0).
\end{eqnarray*}
By the last equality we conclude either that $(\Lambda,E) -
i_{\ast}(\Lambda_N,E_N) = (0,0) \Leftrightarrow (\Lambda,E) =
i_{\ast}(\Lambda_N,E_N)$, i.e. $i : (N, \Lambda_N,E_N) \to
(M,\Lambda,E)$ is a Jacobi map, or that $\Lambda =
i_{\ast}\Lambda_N + \sum_{j=1}^k X_j\wedge Y_j$ and $E =
i_{\ast}E_N +X$, where $X_j, X \in \Gamma(\pi(D))$, $Y_j \in
\Gamma(TM)$, $j=1,\ldots,k$, are convenient vector fields such
that $[\Lambda,\Lambda]=-2E\wedge \Lambda$ and $[E,\Lambda]=0$.

\vspace{3mm} \noindent {\bf Particular cases}
\begin{enumerate}
\item[{\bf a)}] When $D=\{(0,0)\}$, then $D^{\bot}=T^{\ast}M\times
\R$, and they verify the assumptions of Theorem
\ref{th-reduction}. Condition $D_0=D\cap (TN\times \R)=\{(0,0)\}$
is automatically satisfied and the reducibility condition
(\ref{eq4}) takes the form
$$
\Lambda^{\#}(T^{\ast}M)\subseteq TN\;\;\mathrm{on}\;\;N
\hspace{3mm}\mathrm{and}\hspace{3mm} E\vert_N \in \Gamma(TN),
$$
which is exactly the condition given in \cite{dlm} and
\cite{marle} for the submanifolds $N$ of $(M,\Lambda,E)$ of the
first kind.

\item[{\bf b)}] When $D = (\Lambda,E)^{\#}((TN \times
\R)^{\bot})$, we have that $D$ has only sections of type $(X,0)$
if and only if $E\vert_N \in \Gamma(TN)$ and $D_0 = D\cap
(TN\times \R) = \{(0,0)\}$ if and only if $TN\cap
\Lambda^{\#}(TN^{\bot})=\{0\}$. Thus, under the assumptions
\begin{equation}\label{44}
TN\cap \Lambda^{\#}(TN^{\bot})=\{0\} \;\;\mathrm{on}\;\;N
\hspace{3mm} \mathrm{and}\hspace{3mm} E\vert_N \in \Gamma(TN),
\end{equation}
by a simple calculation we show that $D =
\Lambda^{\#}(TN^{\bot})\times \{0\}$ is a Lie subalgebroid of
$(TM\times \R,[\,,\,],\pi)$ if and only if $\Lambda$ belongs to
the ideal generated by the space of smooth sections of $TN$. Also,
since $\Lambda^{\#}((\Lambda^{\#}(TN^{\bot}))^{\bot})\subseteq TN$
and $E\vert_N \in \Gamma (TN)$, it is easy to prove that $D^{\bot}
= (\Lambda^{\#}(TN^{\bot}))^{\bot}\times \R$ is a Lie subalgebroid
of $(T^{\ast}M\times \R, [\,,\,]_{(\Lambda,E)}, \pi\circ
(\Lambda,E)^{\#})$.

\hspace{5mm}Consequently, if (\ref{44}) holds and $\Lambda$
belongs to the ideal generated by the space of smooth sections of
$TN$, then we have that the requirements of Theorem
\ref{th-reduction} as the reducibility condition (\ref{43}) are
verified, therefore $(\Lambda,E)$ induces a Jacobi structure on
$N$. We note that conditions (\ref{44}) are exactly those given in
\cite{i1}.
\end{enumerate}

\vspace{3mm} \noindent {\bf 2. Reduction of Jacobi manifolds with
symmetry :} Let $(M,\Lambda,E)$ be a Jacobi manifold, $G$ a
connected Lie group acting on $M$ by a Jacobi action,
$\mathcal{G}$ the Lie algebra of $G$, $\mathcal{G}^{\ast}$ the
dual space of $\mathcal{G}$ and $J : M \to \mathcal{G}^{\ast}$ an
$Ad^{\ast}$-equivariant moment map for the considering action. Let
$D$ be the vector subbundle of $TM\times \R$ formed by the pairs
$(X_M,0)$, where $X_M$ is the fundamental vector field on $M$
associated to an element $X\in \mathcal{G}$, and $D^{\bot}$ its
conormal bundle which is $D^{\bot}=\{X_M \in TM \,/\,X\in
\mathcal{G}\}^{\bot}\times \R$. It is easy to check that $D$ and
$D^{\bot}$ are Lie subalgebroids of $(TM\times \R,[\,,\,],\pi)$
and $(T^{\ast}M\times \R, [\,,\,]_{(\Lambda,E)}, \pi\circ
(\Lambda,E)^{\#})$, respectively. (For $D^{\bot}$, we take into
account that the action of $G$ on $M$ is a Jacobi action, thus,
for any fundamental vector field $X_M$ on $M$,
$\mathcal{L}_{X_M}\Lambda = 0$ and $\mathcal{L}_{X_M}E = 0$.)
Consequently, $L=D\oplus D^{\bot}$ is a Dirac subbundle of
$((TM\times \R)\oplus (T^{\ast}M\times \R), (0,1)+(-E,0))$. We
suppose that $0$ is a weakly regular value of the moment map $J$.
Hence, $N=J^{-1}(0)$ is a submanifold of $M$ and $D_0=D\cap
(TN\times \R)=\{(X_M,0)\,/\, X\in \mathcal{G}_0\}$, where
$\mathcal{G}_0$ is the Lie algebra of the isotropy subgroup $G_0$
of $0$. Also, we suppose that $\pi(D)$ and $\pi(D_0)$ define,
respectively, a simple foliation $\mathcal{F}$ of $M$ and a simple
foliation $\mathcal{F}_0$ of $N$. Since,
$(\Lambda,E)^{\#}(D^{\bot})\subseteq TN\times \R +D$ holds on $N$,
from the Reduction Theorem \ref{th-reduction} we get that
$(\Lambda,E)$ induces a Jacobi structure on $N/\mathcal{F}_0$. For
more details, see \cite{j3}, \cite{mk} and \cite{ib}.

\vspace{7mm}\noindent {\bf Acknowledgments}

\vspace{3mm}\noindent The authors thank David Iglesias for his
important suggestions. Research partially supported by
GRICES/French Embassy (Project 502 B2) and CMUC-FCT.

\biblio
\bibliographystyle{plain}

\begin{thebibliography}{35}

\bibitem{acw}{A. Cannas da Silva and A. Weinstein, {\it Geometric Models for Noncommutative Algebras},
University of California, Berkeley Mathematics Lecture Notes 10 -
AMS, Providence, 1999.}

\bibitem{cw}{T. Courant and A. Weinstein, {\it Beyond Poisson
structures}, in S\'eminaire Sud-Rhodanien de G\'eometrie, Travaux
en cours 27, Hermann, Paris 1988, pp. 39-49.}

\bibitem{c}{T. Courant, {\it Dirac manifolds}, Trans.
Amer. Math. Soc. 319 (1990), 631-661.}

\bibitem{dlm}{P. Dazord, A. Lichnerowicz, C.-M. Marle, {\it Structure locale des vari\'et\'es de Jacobi}, J. Math. Pures
Appl. 70 (1991) 101-152.}

\bibitem{dr}{V. Drinfeld, {\it Quantum groups}, in Proceedings of the International
Congres of Mathematicians, Berkeley, AMS 1986, pp. 798-820.}

\bibitem{df}{J.-P. Dufour, {\it Normal forms for Lie algebroids}, in Lie Algebroids, Banach Center Publications, Vol. 54,
Warszawa 2001, pp. 35-41.}

\bibitem{fr}{R.L. Fernandes, {\it Lie Algebroids, Holonomy and Characteristic Classes}, Adv.
Math. 170 (2002) 119-179.}

\bibitem{gr}{J. Grabowski, {\it Abstract Jacobi and Poisson structures}, J. Geom. Phys. 9 (1992) 45-73.}

\bibitem{gm1}{J. Grabowski and G. Marmo, {\it Jacobi structures revisited}, J. Phys. A : Math. Gen.
34 (2001) 10975-10990.}

\bibitem{gm2}{J. Grabowski and G. Marmo, {\it The graded Jacobi algebras and (co)homology}, J. Phys. A : Math. Gen.
36 (2003) 161-181.}

\bibitem{ib}{A. Ibort, M. de Leon and G. Marmo, {\it Reduction of Jacobi manifolds},
J. Phys. A : Math. Gen. 30 (1997) 2783-2798.}

\bibitem{im1}{D. Iglesias and J.C. Marrero, {\it Generalized Lie bialgebroids and Jacobi structures},
J. Geom. Phys. 40 (2001) 176-200.}

\bibitem{im3}{D. Iglesias, B. Lopez, J.C. Marrero and E. Padr\'on, {\it Triangular generalized Lie bialgebroids :
Homology and cohomology theories}, in Lie Algebroids, Banach
Center Publications, Vol. 54, Warszawa 2001, pp. 111-133.}

\bibitem{i1}{D. Iglesias Ponte, {\it $\mathcal{E}^1(M)$-Dirac structures and Jacobi structures}, in
Differential geometry and its applications, Proc. Conf. Opava
2001, Silesian Univ. Opava, (Opava, 2001), pp. 275-283.}


\bibitem{im2}{D. Iglesias and J.C. Marrero, {\it Lie algebroid foliations and $\mathcal{E}^1(M)$-Dirac structures},
J. Phys. A : Math. Gen. 35 (2002) 4085-4104.}

\bibitem{krb}{Y. Kerbrat and Z. Souici-Benhammadi, {\it Vari\'et\'es de Jacobi
et groupo\"{\i}des de contact}, C. R. Acad. Sci. Paris, S\'erie I,
317 (1993) 81-86.}

\bibitem{kr}{A. Kirillov, {\it Local Lie algebras}, Russian Math. Surveys 31 (1976) 55-75.}

\bibitem{ks1}{Y. Kosmann-Schwarzbach, {\it Exact Gerstenhaber algebras and Lie bialgebroids}, Acta Appl. Math 41
(1995) 153-165.}

\bibitem{kz}{J.-L. Koszul, {\it Crochet de Schouten-Nijenhuis et cohomologie},
in \'Elie Cartan et les Math\'ematiques d'aujourd'hui,
Ast\'erisque, Num\'ero Hors S\'erie (1985) 257-271.}

\bibitem{lch}{A. Lichnerowicz, {\it Les vari\'et\'es de Jacobi et leurs alg\`ebres de Lie associ\'ees},
J. Math. pures et appl. 57 (1978) 453-488.}

\bibitem{lwx1}{Z.-J. Liu, A. Weinstein, P. Xu, {\it Manin triples for Lie bialgebroids},
J. Diff. Geom. 45 (1997) 547-574.}

\bibitem{lwx2}{Z.-J. Liu, A. Weinstein, P. Xu, {\it Dirac Structures and Poisson Homogeneous Spaces},
Commun. Math. Phys. 192 (1998) 121-144.}

\bibitem{liu}{Z.-J. Liu, {\it Some remarks on Dirac structures and Poisson reductions}, in
Poisson Geometry, Banach Center Publications, Vol. 51, Warszawa
2000, pp. 165-173.}

\bibitem{mck}{K. Mackenzie, {\it Lie groupoids and Lie algebroids in differential geometry}, London Math. Soc.
Lecture notes series 124, Cambridge University Press, Cambridge
1987.}

\bibitem{mx}{K. Mackenzie and P. Xu, {\it Lie bialgebroids and Poisson groupoids},
Duke Math. J. 73 (1994) 415-452.}

\bibitem{marle}{Ch.-M. Marle, {\it On submanifolds and quotients of Poisson and Jacobi manifolds},
in Poisson Geome-try, Banach center publications, Vol. 51,
Warszawa 2000, pp. 197-209.}

\bibitem{mrl}{Ch.-M. Marle, {\it Differential calculus on a Lie algebroid and Poisson manifolds},
The J.A. Pereira da Silva Birthday Schrift, Textos de Matem\'atica
32, Departamento de Matem\'atica da Universidade de Coimbra,
Portugal (2002) pp. 83-149.
(http://www.math.jussieu.fr/$\sim$marle/)}

\bibitem{mr}{J. Marsden and T. Ratiu, {\it Reduction of Poisson manifolds},
Lett. Math. Physics 11 (1986) 161-169.}

\bibitem{mk}{K. Mikami, {\it Reduction of local Lie algebra structures},
Proc. Amer. Math. Soc. 105 (1989) 686-691.}


\bibitem{j2}{J.M. Nunes da Costa, {\it R\'eduction des vari\'et\'es de Jacobi},
C.R.A.S. Paris 308 S\'erie I (1989) 101-103.}

\bibitem{j3}{J.M. Nunes da Costa, {\it Une g\'en\'eralisation, pour les vari\'et\'es de Jacobi,
du th\'eor\`eme de Marsden-Weinstein}, C. R. Acad. Sci. Paris,
S\'erie I, 310 (1990) 411-414.}

\bibitem{j1}{J.M. Nunes da Costa, {\it Compatible Jacobi manifolds : geometry and reduction},
J. Phys. A : Math. Gen. 31 (1998) 1025-1033.}

\bibitem{jj}{J.M. Nunes da Costa and J. Clemente-Gallardo, {\it Dirac structures for generalized Lie
bialgebroids}, J. Phys. A : Math. Gen. 37 (2004) 2671-2692.}

\bibitem{ro}{D. Roytenberg, {\it Courant algebroids, derived brackets and even symplectic supermanifolds},
Ph.D. Thesis, University of California, Berkeley 1999. (arXiv:
Math.DG/9910078.)}

\bibitem{wd}{A. Wade, {\it Conformal Dirac structures}, Lett. Math. Phys. 53 (2000) 331-348.}

\end{thebibliography}
\end{document}